\documentclass[10pt,a4paper,twoside]{amsart}
\usepackage{amssymb,amsmath,amsthm,enumerate,mathrsfs}

\usepackage[]{hyperref}
\hypersetup{
     pdfnewwindow=true,      
     colorlinks=true,       
     linkcolor=black,          
     citecolor=black,        
     filecolor=black,      
     urlcolor=black           
 }

\usepackage[pdftex]{graphics}
\usepackage{graphicx} 

\usepackage[all]{xy}
\newtheorem{thm}{Theorem}[section]
\newtheorem{cor}[thm]{Corollary}
\newtheorem{prop}[thm]{Proposition}
\newtheorem{lemma}[thm]{Lemma}

\newcommand{\bthm}{\begin{thm}}   
\newcommand{\ethm}{\end{thm}}
\newtheorem{pro}[thm]{Proposition} 
\newcommand{\bpro}{\begin{pro}}
\newcommand{\epro}{\end{pro}}
\newtheorem{lem}[thm]{Lemma}
\newcommand{\ble}{\begin{lem}}
\newcommand{\ele}{\end{lem}}
\newcommand{\bco}{\begin{cor}}
\newcommand{\eco}{\end{cor}}
\newcommand{\bprf}{\begin{proof}}
\newcommand{\eprf}{\end{proof}}

\theoremstyle{definition}
\newtheorem{de}[thm]{Definition}
\newcommand{\bde}{\begin{de}}
\newcommand{\ede}{\end{de}}
\newtheorem{ex}[thm]{Example}
\newcommand{\bex}{\begin{ex}}
\newcommand{\eex}{\end{ex}}

\theoremstyle{remark}
\newtheorem{rem}[thm]{Remark}
\newcommand{\bre}{\begin{rem}}
\newcommand{\ere}{\end{rem}}
\newtheorem*{no}{Notes}
\newcommand{\bno}{\begin{no}}
\newcommand{\eno}{\end{no}}

\renewcommand\labelenumi{\theenumi)}

\newcommand{\bb}{\mathbb} 

\newcommand{\R}{\bb{R}}
\newcommand{\Z}{\bb{Z}}
\newcommand{\N}{\bb{N}}

\newcommand{\beq}{\begin{eqnarray*}}
\newcommand{\eeq}{\end{eqnarray*}}
\DeclareMathOperator{\im}{im}
\DeclareMathOperator{\Lip}{Lip}
\DeclareMathOperator{\spt}{spt}
\DeclareMathOperator{\diam}{diam}
\newcommand{\defl}{\mathrel{\mathop:}=}
\newcommand{\defr}{=\mathrel{\mathop:}\,}

\DeclareMathOperator{\id}{id}

\numberwithin{equation}{section}

\begin{document}

\title{Singular (Lipschitz) homology and homology of integral currents}

\author{Christian Riedweg}
\address{Department of mathematics\\ 
ETH Z\"urich\\ 
R\"amistrasse 101\\
8092 Z\"urich\\
Switzerland}
\email{\href{mailto:riedweg@math.ethz.ch}{riedweg@math.ethz.ch}}
\thanks{The first author was partially supported by the Swiss National Science Foundation.
}

\author{Daniel Sch\"appi}
\address{Department of Mathematics\\
University of Chicago\\
5734 S University Avenue\\
Chicago, IL 60637\\
USA}
\email{\href{mailto:schaeppi@math.uchicago.edu}{schaeppi@math.uchicago.edu}}

\keywords{Singular (Lipschitz) homology, homology of compactly supported integral currents}


\begin{abstract}
We compare the homology groups $H_n ^{IC}(X)$ of the chain complex of integral currents with compact support of a metric space $X$ with the singular Lipschitz homology $H^L _n (X)$ and with ordinary singular homology. If $X$ satisfies certain cone inequalities all these homology theories coincide. On the other hand, for the Hawaiian Earring the homology of integral currents differs from the singular Lipschitz homology and it differs also from the classical singular homology $H_n(X)$.
\end{abstract}

\maketitle

\section{Introduction and results} \label{intro}
Building on the theory of metric currents (see \cite{ak00}), Stefan Wenger introduced the homology of integral metric currents with compact support (\cite{wenger06}) in complete quasiconvex metric spaces that admit certain cone inequalities. As the axioms of Eilenberg--Steenrod for a homology theory are satisfied, this homology is isomorphic to the singular (Lipschitz) homology on finite  CW-complexes.

Here, we compare this  homology of metric currents with the singular homology and the singular Lipschitz homology. If locally cycles can be filled by a diameter-controlled chain, all three theories are identical. In the special case of the Hawaiian Earring  we show that the later homology theories do not coincide with the first one.

\bde
A subset $S$ of  a metric space $X$ admits a strong $\gamma$-Lipschitz contraction, $\gamma >0$, if there is a map $\phi:[0,1] \times S \to X$  and $x_0 \in X$ such that $\phi(1, \cdot) = \id_S$, $\phi(0, s) = x_0,\, \forall s \in S,$ and
$$d(\phi(t,s), \phi(t',s')) \le \gamma \diam(S) |t-t'| +\gamma d(s,s').$$

The space $X$ admits locally strong Lipschitz contractions if for all $x \in X$ there is $r_x >0 $ and $\gamma_x >0$ such that every subset $U \subset B_{r_x}(x)$ admits a strong $\gamma_x$-Lipschitz contraction.
\ede

\bre
What we call a strong $\gamma$-Lipschitz contraction is also called simply a {\it $\gamma$-Lipschitz contraction} (e.g.\ in \cite{wenger06}). However we want to emphasize that the contracting map 
has this diameter bound for the first entry.
\ere

For example complete CAT($\kappa$)-spaces admit locally strong $\gamma$-Lipschitz contractions; we come back to this in \autoref{cone}.

We denote the group of $k$-dimensional integral currents in $X$ with compact support by ${\bf I}_k ^c (X)$. We let $C_k^L(X)$ be the group of singular Lipschitz $k$-chains, and we let $C_k(X)$ be the group of singular $k$-chains, both with coefficients in $\Z$. For $c \in C_k ^L(X)$ the current induced by $c$ is denoted by $[c]$. The boundary of a current $T$ is $\partial T$; the boundary of a chain $c$ is $bc$.
For a (Lipschitz) chain $c$, we denote by $sd^m (c)$ the $m$-th barycentric subdivision of $c$. For $A \subset X$, $B_r(A) \defl \{y \in X; \, d(A,y) < r \}$ is the open $r$-neighbourhood of $A$. The support of the current $T$ is denoted by $\spt(T)$, and we write $\im(c)$ for the union of the images of the simplices of a chain. Definitions of these concepts can be found in \autoref{defs}.

We write $H_n ^{IC}(X)$ for the homology of the chain complex of integral currents with compact support of a metric space $X$,  $H^L _n (X)$ for the singular Lipschitz homology and $H_n(X)$ for ordinary singular homology.

\bthm     \label{T1}
Let $X$ be a complete metric space. For $T \in {\bf I}_0 (X)$ there exists $c \in C_0 ^L(X)$ with $T = [c].$

Suppose the complete metric space $X$ admits  locally strong Lipschitz contractions, let $\epsilon >0$ and $n \in \N$. Given   $T \in {\bf I}_n ^c(X)$  with $\partial T = [c]$ for $c \in C_{n-1}^L(X)$ and $bc=0$, there exist $N \in \N$, $T_1, \ldots, T_N \in  {\bf I}_n ^c(X)$, $c_1, \dots, c_N \in C_n^L(X)$ and $V_1, \ldots, V_N \in {\bf I}_{n+1} ^c (X)$ such that
\begin{enumerate}
\item  \label{T1i}
$\sum T_i = T$ and $\partial V_i = T_i - [c_i]$ 
\item
$\exists m \in \N \cup \{0\}$: $b(\sum c_i) = sd^m (c).$
\item
$\spt(V_i) \cup \im(c_i) \subset B_\epsilon (\spt (T)\cup \im(c))$ and $\diam(\spt(V_i)) < \epsilon.$
\end{enumerate}
\ethm

Note that \ref{T1i}) implies that $T$ and $\sum [c_i]$ are homologous.

\bco \label{C4}
If $X$ and the closed subset $A \subset X$ admit  locally strong Lipschitz contractions then the homology of integral Lipschitz chains is isomorphic to the homology of integral currents with compact support:
$$H_k  ^{L}(X, A) \	\cong  H_k ^{IC} (X,A).$$
The isomorphism is induced by the map $[\,]:C_*^L(X) \to  {\bf I}_* ^c(X)$.
\eco

\bthm     \label{Tsing}
Let $X$ be a metric space. Then  $C_0 ^L(X) = C_0(X).$

Suppose that the  metric space $X$ admits  locally strong Lipschitz contractions, let $\epsilon >0$ and $n \in \N$. Given   $c \in C_n(X)$  with $bc  \in C_{n-1}^L(X)$, there exist $N \in \N$, $c_1 , \ldots, c_N \in C_n (X)$, $c^L _1, \dots, c^L _N \in C_n^L(X)$ and $\bar c_1, \ldots, \bar c_N \in C_{n+1} (X)$ such that
\begin{enumerate}
\item $\exists m \in \N \cup \{0\}$: $\sum c_i = sd^m(c)$ and $b \bar c_i = c_i - c^L _i$ 
\item
$\im(\bar c_i)  \subset B_\epsilon (\im(c))$  and $\diam(\im(\bar c_i)) < \epsilon.$
\end{enumerate}

\ethm

\bco \label{C5}
If $X$ and the  subset $A \subset X$ admit  locally strong Lipschitz contractions then the homology of integral Lipschitz chains is isomorphic to the singular  homology:
$$H_k  ^{L} (X, A) \	\cong  H_k (X,A).$$
The isomorphism is induced by the inclusion $C^L_*(X) \subset C_*(X).$
\eco

The above theorems are proved in a more general setting where we make assumptions on local $k$-cycles for $0 \le k \le n$, see \autoref{proofT1}.

The \emph{Hawaiian Earring} $\mathbb{H} \subset \mathbb{R}^2$ is given by the countable union of the circles
\[
L_n=\{x \in \mathbb{R}^2; \left\| x-(1/n,0) \right\|=1/n\}
\]
with radius $1/n$ and centre $(1/n,0)$, $n\in \mathbb{N}$. As metric on $\mathbb{H}$ we set
\[
 d(x,y) \defl \left\{ \begin{matrix} \left\|x-y\right\|, & \text{ if } \exists n \in \N: \, x,y \in L_n,\\
                       \left\|x\right\| + \left\|y\right\|, & \text{ otherwise}
                      \end{matrix} \right.
\]
(we could as well take the length metric of $\mathbb H$). Note that any neighbourhood of $(0,0)$ in $\mathbb{H}$ contains all but finitely many of the $L_n$. Thus $\mathbb{H}$ is not locally contractible, and in particular not a CW-complex. 

We show in \autoref{proofHawaiian} that the maximal divisible subgroup of $H_1^{IC}(\mathbb{H})$ is trivial whereas   the maximal divisible subgroups of $H_1^L(\mathbb{H})$ and $H_1(\mathbb{H})$ are non-trivial. This implies that
\bthm\label{HawaiianEarring}
$H_1^{IC}(\mathbb{H})$ is  isomorphic neither to $H_1^L(\mathbb{H})$ nor to $H_1(\mathbb{H})$.
\ethm

Note that Thierry De Pauw  compares in \cite{DePauw} various homology theories. However, neither metric currents nor singular Lipschitz chains are considered therein.

The paper is organized as follows: The definitions of currents  and singular (Lipschitz) chains are given in \autoref{defs}; there we also recall facts needed later on.  
The proofs of Theorem~\ref{T1} and Theorem~\ref{Tsing} as well as of their corollaries are given in \autoref{proofT1}; in \autoref{proofHawaiian} we prove Theorem~\ref{HawaiianEarring}.

\subsection*{Acknowledgements}{We are grateful to Urs Lang, Roger Z\"ust and Lorenz Reichel for helpful discussions and comments.}

\section{Definitions} \label{defs}

\subsection{Currents} \label{currents} We recall some definitions from \cite{ak00}. Let $X$ be a complete metric space.
By $D^k (X)$ we denote the set of $(k+1)$-tuples $(f,\pi_1, \ldots, \pi_k)$ of Lipschitz functions from $X$ to $\R$ with $f$ bounded.
We write $\Lip(f)$ for the Lipschitz constant of the Lipschitz function $f{:}\, X \to \R$.
\bde[{\cite[2.1]{wenger06}}]
A $k$-dimensional current $T$ on $X$ is a multi-linear functional on $D^k(X)$
satisfying the following properties:
\begin{enumerate}
\item 
If $\pi^j _i$ converges point-wise to $\pi_i$ as $j \to \infty$ and if $\sup_{i,j} \Lip(\pi^j _i) < \infty$ then
$T(f,\pi^j _1, \ldots, \pi ^j _k) \to T(f,\pi _1, \ldots, \pi  _k).$
\item \label{locality}
If $\{x \in X; \, f(x) \neq 0   \}$ is contained in the union $\bigcup_{i=1} ^k B_i$ of Borel sets $B_i$ and if $\pi_i$ is constant on $B_i$ then 
$T(f,\pi _1, \ldots, \pi  _k) =0.$
\item
There exists a finite Borel measure $\mu $ on $X$ such that
\begin{eqnarray} \label{mass}
|T(f,\pi _1, \ldots, \pi _k)| \leq \prod _{i=1} ^k \Lip(\pi _i) \int_X |f|d\mu
\end{eqnarray}
for all $(f, \pi_1, \ldots , \pi_k) \in D^k (X).$
\end{enumerate}
Let ${\bf M}_k(X)$ be the set of $k$-dimensional  currents. 
\ede

The minimal Borel measure $\mu$ with (\ref{mass}) is called mass of $T$ and written as $\|T\|$; set $M(T) \defl \|T\|(X)$. 

The support $\spt(T)$ of $T$ is the closed  set given by
 $$\spt(T) \defl	\{ x \in X; \, \|T\|( B_r(x)) > 0,\, \forall r>0\}.$$

\bre
As is done in \cite{ak00, wenger06}, we will henceforth assume that the cardinality of $X$ is an Ulam number. Then for $T \in {\bf M}_k (X)$ one has 
$\|T\|(X \backslash \spt(T)) =0$ (see \cite{ak00}).
\ere

Let $\chi _C$ be the characteristic function of a set $C$.
For a Borel set $A \subset X$ and $T \in {\bf M}_k(X)$ we define the restriction of $T$ to $A$ as
$$(T \lfloor A) (f, \pi_1 ,\dots, \pi_k) \defl T(\chi_A f, \pi_1 ,\dots, \pi_k),$$
which is in ${\bf M}_k(X)$  (the above expression is well-defined since $T$ extends to a functional whose first argument lies in $L^\infty (X, ||T||),$ see \cite{ak00}). 

The boundary of $T \in {\bf M}_k(X)$ is the functional
$$\partial T (f, \pi_1, \ldots, \pi_{k-1}) \defl T(1, f, \pi_1, \ldots, \pi_{k-1}).$$
We say that $T$ is normal if $\partial T \in {\bf M}_{k-1}(X)$; the space of normal currents is denoted by ${\bf N}_k(X).$ From the locality property for currents (condition \ref{locality})) we get $\partial (\partial T) =0.$

For a Lipschitz map $g: X \to Y$ to a complete metric space $(Y,d)$ and a current $T \in {\bf M}_k(X)$ we define the push-forward of $T$ by $g$ to be
$g_{\#}T (f, \pi_1, \ldots, \pi_k) \defl T(f\circ g, \pi_1 \circ g, \ldots , \pi_k \circ g)$
for $(f, \pi_1, \ldots , \pi_k) \in D^k(Y).$ Then, $\partial(g_{\#}T ) = g_{\#}(\partial T).$

\bde[{\cite[2.3]{wenger06}}]
A  current $T \in  {\bf M}_k(X)$ with $k \geq 1$ is said to be rectifiable if
\begin{enumerate}
\item
 $||T||$ is concentrated on a countably ${\mathcal H}^k$-rectifiable set and
\item
$||T||$ vanishes on ${\mathcal H}^k$-negligible sets.
\end{enumerate}
The current $T$ is called integer rectifiable if, in addition, the following property holds:
\begin{enumerate}\setcounter{enumi}{2}
\item
For any Lipschitz map $\phi: X \to \R^k$ and any open set $U \subset X$ there exists $\theta \in L^1(\R^k, \Z)$ such that
$$\phi_{\#}(T \lfloor  U)(f, \pi_1, \ldots , \pi _k) = \int_{\R^k} \theta \, f \, \det( \frac {\partial \pi_i}{\partial s_j}) d {\mathcal L}^k (s)$$
for all $(f, \pi_1, \ldots, \pi_k) \in D^k(X).$
\end{enumerate}
A $0$-dimensional (integer) rectifiable current is a $T \in {\bf M}_0 (X)$ of the form
\begin{eqnarray}
\label{0Current}
T(f) = \sum _{i=1} ^{\infty} \theta_i f(x_i), 
\end{eqnarray}
$f$ Lipschitz and bounded, for suitable $\theta_i \in \R$ ($\theta_i \in \Z$ respectively, note that in this case the sum has to be finite) and $x_i \in X.$ We write $T = \sum_i \theta_i [x_i]$.
\ede

The space of integral currents, i.e.\ integer rectifiable and normal currents (equivalently defined as integer rectifiable currents with integer rectifiable boundary, see \cite[Theorem 8.6]{ak00}), is denoted by ${\bf I}_k(X)$ and
$${\bf I}_k ^c (X) \defl \{T \in {\bf I}_k(X); \, \spt(T) \mbox{ is compact} \}$$
is the space of integral currents with compact support.
So, $\partial |_{{\bf I}_k(X)}{:}\, {\bf I}_k(X) \to {\bf I}_{k-1}(X)$,  $\partial|_{{\bf I} ^ c_k(X)}{:}\, {\bf I} ^ c_k(X) \to {\bf I} ^ c_{k-1}(X)$ and 
 $g_{\#}T \in {\bf I}^c_k(Y)$  for $T \in {\bf I}^c_k(X)$ and $g$ Lipschitz.

Let $A \subset X$ be a closed subset, we define (see \cite[p. 159]{wenger06})
$${\mathcal Z}_k ^{IC}(X,A) \defl \{T \in {\bf I}_k ^c (X); \partial T \in {\bf I}_{k-1} ^c (A) \}$$
$${\mathcal B}_k ^{IC}(X,A) \defl \{R + \partial S; R \in {\bf I}_k ^c (A), \, S \in {\bf I}_{k+1}^c (X)\}.$$
The homology of integral currents with compact support is 
$$H^{IC} _k (X, A) \defl {\mathcal Z}_k ^{IC}(X,A)  / {\mathcal B}_k ^{IC}(X,A).$$
If $A = \emptyset$ we write $H^{IC} _k (X)$.

For $T \in {\bf N}_k(X)$ and a Lipschitz function $d: X \to \R$ we set
\beq
 \langle T, d, r+\rangle &\defl& \partial(T\lfloor{\{d \leq r\}} ) - (\partial T)\lfloor{ \{d \leq r\}}\\
&=& (\partial T)\lfloor{ \{d > r\}}  - \partial (T\lfloor{\{d > r\}}).
\eeq
If $T \in {\bf I}^c_k(X)$ then,  for almost all $r \in \R$, $\langle T,d,r+\rangle \in {\bf I}^c_{k-1}(X)$. This is  called the slice of $T$ by $d$ at $r$ and we  denote it by
$\langle T, d, r\rangle$.
Note that
$\langle \partial T, d, r+\rangle = - \partial \langle T, d, r+\rangle$ and that 
$$\spt(\langle T, d, r \rangle) \subset d^{-1} \{r\} \cap \spt(T).$$

\subsection{Singular (Lipschitz) chains} \label{chains}
Let 
$$\Delta ^k \defl \{(s_0, \ldots , s_k) \in \R^{k+1}; \, \sum_j s_j =1 \text{ and } 0 \leq s_j, \, \forall j   \} \subset \R^{k+1}$$
 be the standard simplex. Sometimes it is more convenient to consider $\Delta^k$ as a subset of $\R^k$; for this we choose an isometry $\phi: \R^k \to \{s \in \R^{k+1}; \, \sum_j s_j =1\}$ and define $\tilde{\Delta}^k = \phi^{-1}(\Delta ^k).$ Let $X'$ be a metric space (not necessarily complete).

\bde
A singular  $k$-simplex $c$ is a continuous map $c: \Delta^k \to X'$.

A singular  $k$-chain $c$ over an abelian group $G$ is a finite formal sum $c= \sum_{i=1} ^m n_i c_i$, where $n_i \in G$ and $c_i$ are singular simplices.

A singular Lipschitz $k$-chain is a singular $k$-chain $c= \sum_{i=1} ^m n_i c_i$ where all $c_i$ are Lipschitz maps.
\ede

We use only $G= \Z$ and speak of integral (Lipschitz) $k$-chains. We denote by $C_k (X')$ 
 the free abelian group of integral $k$-chains and by $C_k ^L(X')$
the subgroup of integral Lipschitz $k$-chains.

Let $[e_0,e_1,\ldots,e_k]$ be the vertices of the standard $k$-simplex for $k \geq 1$, then the boundary of a singular (Lipschitz)  $k$-simplex $c$ is the (Lipschitz) $(k-1)$-chain 
$$bc \defl \sum_{j=0} ^k (-1)^j c _{|_{[e_0,\ldots,\hat{e}_{j},\ldots e_k] }}$$
 where $\hat{e}_j$ means that $e_j$ is omitted.
By setting $bc = \sum_{i=1} ^m n_i b c_i$ for a $k$-chain $c=\sum_{i=1} ^m n_i c_i$ and defining  $bc=0$ for $0$-chains, we get homomorphisms $b: C_k ^{(L)}(X') \to C_{k-1}^{(L)}(X')$ such that $b(bc)=0$. 
For a  subset $A \subset X'$ let
$${\mathcal Z}_k ^{(L)}(X',A) \defl \{c \in C^{(L)} _k (X'); \, bc \in C^{(L)}_{k-1} (A) \}$$
$${\mathcal B}_k ^{(L)}(X',A) \defl \{c + b \bar c; \, c \in C^{(L)}_k (A), \, \bar c \in C^{(L)}_{k+1} (X') \}$$
and
$$H^{(L)} _k (X', A) \defl {\mathcal Z}_k ^{(L)}(X',A) /{\mathcal B}_k ^{(L)}(X',A).$$
If $A = \emptyset$ we write $H^{(L)} _k (X')$. There is a natural comparison homomorphism $H_n^L(X, A)\rightarrow H_n(X,A)$ which is induced by the inclusions $C^{L}_n(S)\rightarrow C_n(S)$ for subsets $S \subset X$.

For the complete metric space $X$, any integral Lipschitz $k$-chain $c=\sum_{i=1} ^m n_i c_i$ induces an integral current $[c]$ with compact support  defined by
$$[c] (f, \pi_1, \ldots, \pi_k) \defl \sum_{i=1} ^m n_i \int_{\tilde{\Delta} ^k \subset \R ^k} f \circ \bar{c}_i \, \det \left(\frac {\partial (\pi_j \circ \bar{c}_i)}{\partial s_l}\right) d{\mathcal L}^k(s),$$
where  $\bar{c}_i = c_i \circ \phi  $.  The maps  $[\,]: C_k^L(X) \to {\bf I}_k  ^c(X)$  are homomorphisms. By Stokes' theorem we get  
$[bc] = \partial [c]$, so $[\,]$ is a chain map for these complexes. We also refer to the induced maps between the homologies of these complexes as comparison maps.

For $c= \sum_{i=1} ^m n_i c_i \in C_k(X')$ with $n_i \neq 0$  we define
$$\text{im}(c) \defl  \bigcup _{i =1} ^m \text{im} (c_i)$$
(this is well-defined since $C_k(X')$ is free over the $k$-simplices, so the representation of $c$ as sum of simplices is unique).

The push-forward of $c \in C_k(X')$ to the metric space $Y$ by the continuous map $g{:}\, X' \to Y$ is the chain
$g_{\#} c \defl \sum_{i=1} ^m n_i (g \circ c_i) \in C_k (Y).$ 
Again, $b(g_{\#} c) = g_{\#} (bc)$.
If $c \in C^L_k(X')$ and $g: X' \to Y$ is Lipschitz then $ g_{\#} c \in   C_k ^L (Y)$;
if $X$ and $Y$ are  complete and  $g: X \to Y$ is Lipschitz, $[g_{\#} c] = g_{\#} [c]$ for all $c  \in C^L_K(X).$

In order to get smaller simplices (i.e.\ simplices with smaller diameter) we use barycentric subdivision. This standard construction can be found in \cite{hatcher} or \cite{munkres}. We only list the facts that we need later.
Let $m $ be given and $sd^m(c)$ denote the singular (Lipschitz) $k$-chain  resulting from the $m$-th barycentric subdivision of the singular (Lipschitz) $k$-chain $c$, then:
\begin{enumerate}
\item
For  $k \ge 0$ there is a homomorphism $D_{m,X'}{:} \, C_k ^{(L)}(X') \to C_{k+1} ^{(L)}(X')$ such that for each $k$-chain $c$
\begin{eqnarray} \label{sd}
b(D_{m,X'} (c)) + D_{m,X'}(bc) = sd^m (c) - c.
\end{eqnarray}

Furthermore, $D_{m,X'}$ is natural: for $f{:}\, X' \to Y$ Lipschitz, $f_{\#} \circ D_{m,X'} = D_{m,Y} \circ f_{\#}.$
\item
Applying iterated barycentric subdivision, we can get arbitrary small diameter of the  image of the resulting simplices. 
\item
For $c \in C^L _k(X')$ holds
$[sd^m(c)] = [c].$
\item
$b(sd^m(c)) = sd^m(b(c)).$
\end{enumerate}

\subsection{Cone inequalities}
Here we give the definitions of miscellaneous cone inequalities. These inequalities are used in  Theorem~\ref{T2} and Theorem~\ref{Tsing2} which generalize Theorem~\ref{T1} and~\ref{Tsing}. We recall the  construction of $[0,1] \times T$ for a normal current $T$ (from \cite{wenger06} which is a modified version of the one in \cite{ak00}) and something similar for chains.
This allows us to show in Proposition~\ref{contr-cone} that spaces which admit locally strong Lipschitz contractions also admit these cone inequalities.

\bde
Let $k \ge 1$ and let 
$${\mathcal F} \defl \{F: \R \to \R; \, F \text{ is continuous and non-decreasing with } F(0)=0\}.$$
\begin{itemize}
\item
 $X$ admits a local cone inequality for ${\bf I}_k ^c(X)$  if for every $x \in X$ there exists $r_x >0$ and $F_x \in {\mathcal F}$ such that for every $T \in {\bf I}_k ^c(X)$ with $\partial T = 0$ and $\spt(T) \subset B_{r_x}(x)$ there exists a $\bar T \in {\bf I}_{k+1} ^c(X)$ satisfying $\partial \bar T = T$ and $\diam(\spt(\bar T)) \leq F_x(\diam(\spt(T))).$
\item
$X$ admits a local cone  inequality for  $C_k^{(L)}(X)$ if for every $x \in X$ there exists $r_x >0$ and $F_x\in {\mathcal F}$ such that for every $c \in C_k ^{(L)}(X)$ with $bc = 0$ and $\im(c) \subset B_{r_x}(x)$ there exists a $\bar c \in C_{k+1}^{(L)}(X)$ satisfying $b \bar c = c$ and $\diam(\im(\bar c)) \leq F_x(\diam(\im(c))).$
\item
$X$ admits a local cone inequality for $C_0 ^{(L)}(X)$ if for every $x \in X$ there exists $r_x >0$ and $F_x\in {\mathcal F}$ such that for every $c = \sum n_i c_i \in C_0 ^{(L)}(X)$ with $\sum n_i = 0$ and $\im(c)\subset B_{r_x}(x)$ there exists an $\bar c \in C_{1}^{(L)} (X)$ satisfying $b \bar c = c$ and $\diam(\im(\bar c)) \leq F_x(\diam(\im(c))).$
\end{itemize}
\ede
Note that our definition of cone type inequalities for currents is different from the one used in \cite{wenger06}; there one has the  condition that locally there exists a filling with controlled mass whereas we consider only compactly supported currents and our condition is that locally there exists a filling with controlled diameter.

The following lemma is from \cite[Lemma 30.6]{munkres} (therein it is stated for the singular theory in topological spaces, see below). 
\ble \label{LIxchain}
There exists, for each metric space $X$ and each non-negative integer $k$, a homomorphism
$$
K_X : C_k^{(L)}(X) \to C_{k+1} ^{(L)}([0,1] \times X )$$
having the following property:
If $c \in C_k^{(L)}(X)$ is a singular simplex, then
\begin{eqnarray}\label{Ixchain}
bK_X (c) +K_X(bc) = j_\#(c) - i_\#(c). \end{eqnarray}
Here the map $i: X \to [0,1]\times X $ carries $x$ to $(0,x)$; and the map $j:X \to [0,1]\times X$ carries $x$ to $(1,x)$.
\ele
We only adumbrate the proof (carried out in \cite{munkres} for the continuous case) to indicate that this holds for Lipschitz chains too: One wants to look at $[0,1] \times c$ as a chain in $C_{k+1}^{(L)}([0,1]\times X)$, where $c$ is  a singular simplex in $C_k ^{(L)}(X)$. To do this,  one first gives a  decomposition of $[0,1] \times \Delta^k$ into a $(k+1)$-chain consisting of (regular) simplices in $\R^{k+2}$. Then one carries this decomposition over to $[0,1]\times X$ for every simplex $c \in C_k ^{(L)}(X)$ in an intuitive way, producing a chain $\bar c$ in $C_{k+1} ^{(L)}([0,1] \times \im(c) )$. Clearly,  this construction respects the Lipschitz continuity.

We can also define $[0,1] \times T$ for a normal current $T$ (see \cite[p. 146]{wenger06}): Given a function $f: [0,1] \times X \to \R$ we set $f_t(x)\defl f(t,x),$ so $f_t$ is a function from $X$ to $\R$. For $T \in {\bf N}_k(X)$ we define the normal $k$-current $[t]\times T$ on $[0,1]\times X$ by
$([t]\times T)(f, \pi_1, \ldots, \pi_k)  \defl T(f_t, \pi_{1\, t}, \ldots, \pi_{k \, t}).$
\bde[\cite{wenger06}, Definition 3.1]
For a normal current $T \in {\bf N}_k(X)$ the functional $[0,1] \times T$ on $D^{k+1}([0,1] \times X)$ is given by
\begin{equation*}
 \begin{split}
  ([0,1] \times T)&(f, \pi_1, \ldots, \pi_k) \defl \\
  & \sum _{i=1} ^{k+1} (-1)^{i+1} \int\limits _0 ^1 T(f_t \frac {\partial \pi_{i \,t}}{\partial t}, \pi_{1\,t}, \ldots,\pi_{i-1 \, t}, \pi_{i+1 \, t }, \ldots ,\pi_{k+1 \, t}) dt
 \end{split}
\end{equation*}
for $(f, \pi_1, \ldots, \pi_{k+1}) \in D^{k+1}([0,1] \times X).$
\ede

\bpro[\cite{wenger06}, Theorem 3.2] \label{Ixcurrent}
For $T \in {\bf N}_k (X), $ $k \geq 1$, with bounded support the functional $[0,1] \times T$ is a $(k+1)$-dimensional normal current on $[0,1] \times X$ with boundary
$\partial([0,1] \times T) = [1] \times T - [0] \times T - [0,1] \times \partial T.$
Moreover, if $T \in {\bf I}_k(X)$ then $[0,1] \times T \in {\bf I}_{k+1}([0,1] \times X).$
\epro

\section{Proof of Theorem 1.3, Theorem 1.5 and their corollaries} \label{proofT1}

\subsection{Lipschitz contractions and cone inequalities} \label{cone}
Examples of spaces that admit locally strong Lipschitz contractions are Banach spaces or CAT($\kappa$)-spaces for $\kappa \in \R$. This is discussed on pp.\ 146 and 147 in \cite{wenger06}; as well from there follows the last statement of the following proposition.

\bpro \label{contr-cone}
Let $X$ be a  metric space that admits  locally strong Lipschitz contractions. Then $X$ admits local cone inequalities for $C_j^{L}(X)$ and $C_j(X)$ for all  $j \ge 0.$

If in addition $X$ is complete then it also admits local cone inequalities for ${\bf I}_k ^c(X)$,  $k \ge 1$.
\epro

\bprf
Denote by $x_0 ^n$ the constant $n$-simplex with image $x_0 \in X$. Note that
$bx_0 ^n = 0 $ if $n$ is  odd or zero and $bx_0 ^n = x_0 ^{n-1}$ if $n$ is even.

Let $x \in X$ with $r_x >0$ and  $\gamma_x>0$ as in the definition of  locally strong Lipschitz contractions. Note that for a $\gamma_x$-contraction $\phi$ of a set $S$ we have 
$$\diam(\phi([0,1] \times S)) \le 2 \gamma_x \diam(S).$$
Now we use the homomorphism $K\defl K_X: C_k^{(L)}(X) \to C_{k+1} ^{(L)}([0,1] \times X )$ from Lemma~\ref{LIxchain}.
Let $c \in C_k^{(L)}(X)$ with  $\im(c) \subset B_{r_x}(x)$. Let $\phi: [0,1] \times \im(c) \to X$ be a $\gamma_x$ contraction.
The push-forward of $K(c)$ by $\phi$ is clearly in $ C^{(L)} _{k+1}(X)$. 
If $bc =0$ we have $K(bc) =0$. In this case, for $c = \sum n_i c_i$, we have by (\ref{Ixchain})
$$b(\phi_{\#} K(c)) =c-\sum n_i x_0 ^k .$$
Since $bc = 0$ we  have, for $k$ even and positiv, $\sum n_i = 0$; i.e.\ $\bar c \defl  \phi_{\#}K(c)$ is a  filling of $c$ in this case.  If $k$ is odd, a  filling of $c$ is given by 
\[
 \bar c \defl  \phi_{\#} K(c) + \sum n_i x_0 ^{k+1} 
\]
Note that, if $c$ is a Lipschitz chain, so is $\bar c$. If $k =0$, $c = \sum n_i x_i$ with $\sum n_i =0$ a Lipschitz filling is given by the $1$-chain $\bar c (t) \defl \sum n_i \phi(t, x_i).$
Furthermore,
\[
 \diam(\im(\bar c) ) = \diam(\im(\phi_{\#}K(c))) \le 2 \gamma_x \diam(\im(c)).
\]

Now let $X$ be  complete.
If $T \in {\bf I}_k(X)$ has $\spt(T) \subset B_{r_x}(x)$ and $\partial T =0$, we get a filling $\bar T \defl \phi_{\#} ([0,1] \times T)$ in ${\bf I}_k (X)$ (by Proposition~\ref{Ixcurrent}) with $\spt(\bar{T}) \subset \im(\phi)$, i.e.\ $\diam(\spt(\bar{T})) \leq 2 \gamma \diam(\spt(T))$. If in addition $ T \in {\bf I}_k ^ c(X)$ then $[0, 1] \times T \in {\bf I}_{k+1} ^c ([0,1]\times X)$ and therefore  $\spt(\bar T)$ is compact.

Concluding we see that if  $X$  admits locally strong Lipschitz contractions then $X$ admits local cone inequalities for ${\bf I}_k ^ c(X)$, $k \ge 1$, and $C^{(L)}_j(X)$, $j \ge 0$, with $F_x(t)=2\gamma_x t$. 
\eprf

\subsection{Proof of Theorem~\ref{T1}}
Theorem~\ref{T1} follows by Proposition~\ref{contr-cone} from the more general fact stated below.

\bthm     \label{T2}
Let $X$ be a complete metric space. Then  for  $T \in {\bf I}_0 (X)$ there exists $c \in C_0(X)$ with $T = [c].$

Suppose the complete metric space $X$ admits  local cone inequalities for ${\bf I}_j ^c (X)$ and $C_k^L(X)$ for $j=1, \dots, n$ and $k = 0, \ldots , n-1$ and let $\epsilon >0$. Given   $T \in {\bf I}_n ^c(X)$  with $\partial T = [c]$ for $c \in C_{n-1}^L(X)$ and $bc=0$, there exist $N \in \N$, $T_1, \ldots, T_N \in  {\bf I}_n ^c(X)$, $c_1, \dots, c_N \in C_n^L(X)$ and $V_1, \ldots, V_N \in {\bf I}_{n+1} ^c (X)$ such that
\begin{enumerate}
\item  \label{T2i}
$\sum T_i = T$ and $\partial V_i = T_i - [c_i]$ 
\item \label{T2ii}
$\exists m \in \N _0$: $b(\sum c_i) = sd^m (c).$
\item \label{T2iii}
$\spt(V_i) \cup \im(c_i) \subset B_\epsilon (\spt (T)\cup \im(c))$ and $\diam(\spt(V_i)) < \epsilon.$
\end{enumerate}

\ethm

\bre \label{partofchain}
By a {\it part of a chain} $c = \sum_{i=1}^m n_i c_i$ we mean a chain $c'= \sum_{i=1}^m n'_i c_i$ such that $n_i = n'_i$ or $n'_i = 0$.
Let $c \in C_n ^{(L)} (X)$, $n  \geq 1$ and  $U \subset X$. For $\epsilon >0 $ there exist  a chain $\bar c \in C_n ^{(L)}( B_{\epsilon}(U))$ that is a part of $sd^m(c)$ for some  $m \in \N_0$ and such that $\im(sd^m (c)-\bar c) \subset X\backslash U$. To see this, it is enough to consider a singular $n$-simplex $c$. 
Let $m \in \N$ be such that $sd^m(c)= \sum m_i' c_i'$ with $\diam(c_i') < \epsilon$ for all $i$. Set $\bar c \defl \sum m_i c_i' \in C_n ^{(L)}( B_{\epsilon}(U))$, where $m _i = m_i '$ if  $\im(c_1') \cap U \neq \emptyset$ and $m _i =0$ otherwise.

This fact will be used in the proof to construct a simplicial boundary close to a slice of a simplicial cycle with a filling of the difference: 
Let $bc=0$ and $U= \{d \le r\}$ for a Lipschitz function $d: X \to \R$ and $r \in \R.$ Then, since $[c]\lfloor{U} = [\bar c] - [\bar c]\lfloor (X\backslash U) $, we have
$\langle [c], d, r+ \rangle =  \partial ([c]\lfloor U ) = [b\bar c] -\partial([\bar c] \lfloor (X\backslash U)).$
Thus, $[\bar c] \lfloor (X\backslash U)$ has boundary
\begin{equation} \label{slicetoboundary}
 \partial ([\bar c] \lfloor (X\backslash U)) = [b\bar c] - \langle [c], d, r+ \rangle.
\end{equation}
\ere

\bprf[Proof of Theorem \ref{T2}.]
We argue inductively on the dimension of the current; for the induction step we use another induction on the number of balls needed to cover  $\spt(T) \cup \im(c)$. Let $n=0.$ The decomposition (\ref{0Current}) gives the desired equality for integral currents, i.e.\ $T = \sum_{i=1} ^m n_i [q_i]$ and  $c \defl \sum_{i=1} ^m n_i q_i$.

Now let $n>0:$  For $x \in X$, let $r_x > 0$ denote the minimum of all radii in the (finitely many) cone inequalities for $x$ and $F_x$ the maximum function of all those diameter functions for $x$. 
We can assume that $F_x(r) \ge r.$ 
Let
$0 < R_x <r_x$ be such that 
$$\begin{array}c
2( R_x +  F_x(2R_x + 2 F_x(2R_x))) < \epsilon.
\end{array}$$
Cover  $\spt(T) \cup \im(c)$ by balls of radius $R_x$  and centers in $\spt(T) \cup \im(c)$. We get a finite subcover with centers $x_1, \ldots , x_M$; set $R_i \defl R_{x_i}$ and $F_i \defl F_{x_i}.$ 
We show that there are $V_i,$ $c_i$ and $T_i$, $i =1, \dots , M$, with properties \ref{T2i}) and \ref{T2ii}) of Theorem~\ref{T2} that fulfill moreover
\begin{equation}\label{T2iii'}
\spt(V_i) \cup \im(c_i) \subset B_{\frac \epsilon 2} (x_i);
\end{equation}
this clearly implies \ref{T2iii}).

If $M= 1$  the cone inequality for $C_{n-1} ^L(X)$ gives a filling $c_1 \in C_n ^L(X)$ of $c$ (note that for $n=1$, $c = \sum n_i q_i$ necessarily $\sum n_i = T(1,1) =0)$. We have now $\im(c_1) \subset B_{R_1 + F_1(2R_1)}(x_1)$. Set  $T_1 \defl T$; now the cone inequality for ${\bf I}_n ^c(X)$ gives a filling $V_1 \in {\bf I}_{n+1} ^c(X)$ of $T - [c_1]$ with support in $B_{\frac{ \epsilon}{ 2}}(x_1)$.

If $ M >1$ let  $0 < \bar R < R_M$ be such that 
$\spt(T) \cup \im(c) \subset B_{\bar R}(x_M) \cup \bigcup_{i=1} ^{M-1} B_{R_i} (x_i).$ Set $\alpha \defl \frac{R_M - \bar R}{4}$ and choose $r \in (\bar R + \alpha, \bar R + 3\alpha)$; note that $\bar R < r - \alpha < r +\alpha < R_M.$
We slice $T$  by $d(x) \defl d(x_M, x)$ at $r$. 
We can assume that $\langle T, d, r \rangle  \in {\bf I}_{n-1} ^c(X)$ (and hence
$\langle \partial T, d, r \rangle  \in {\bf I}_{n-2} ^c(X)$ for $n>1$) and that  $\left\| \partial T\right\| (d^{-1}(r)) =0$ for $n=1.$

For $n=1$ set directly $S \defl T\lfloor {\{ d >r\}}$; this is an integer rectifiable $1$-current with compact support. Since  $\partial S = (\partial T)\lfloor {\{ d > r\}} - \langle  T, d, r \rangle \in {\bf I}_0 (X)$ we have $S \in {\bf I}_1 ^c(X).$ As above we see that $\partial S = \sum n_i [q_i]$ with $\sum n_i =0.$  By induction there are 
$c_1, \ldots, c_{M-1},$ $T_1, \ldots, T_{M-1}$ and $V_1, \ldots, V_{M-1}$ with $b \sum c_i = \sum n_i q_i$, property \ref{T2i}) of Theorem~\ref{T2} for $S$ instead of $T$ and with (\ref{T2iii'}). 
Now, $T_{M} \defl T-S$ has $\spt(T_{M}) \subset B_{R_M}(x_M)$ and boundary $\partial (T_{M})= \partial T - \partial S \defr  \sum n'_j [p_j]$; again, $\sum n_j ' =0$. So there is a $ c_{M} \in C_1 ^L (B_{R_M+F_M(2R_M)}(x_M))$ and a filling $V_{M}$ of $T_{M} - [c_{M}]$ with
 $\spt(V_{M}) \sup \im (c_M) \subset 
B_{R_M + F_M(2(R_M + F_M(2R_M)))}(x_M) \subset
B_{\frac \epsilon 2}(x_M)$, proving the case $n=1$.

If $n >1$, choose $\epsilon' \defl \alpha /2$. We can apply  remark~\ref{partofchain} for $c$, $\epsilon'$ and the slice by $d$ at $r$ to get $\bar c \in C_{n-1}^L(B_{r+\epsilon'}(x_M))$ and $m_1 \in \N_0$ with $\im(sd^{m_1}(c) -\bar c) \subset X\backslash B_r (x_M)$.

Now, $T' \defl  \langle T, d, r\rangle - [\bar c]\lfloor {\{d > r\}}  $ has by (\ref{slicetoboundary})
$$\partial T' = - \langle \partial T, d,r\rangle - \partial  ([\bar c]\lfloor {\{d > r\}}) = - [b\bar c],$$
in particular, $T' \in {\bf I}_{n-1} ^c (X)$. By induction assumption for $T'$ and $\epsilon'$ and $c' \defl - b\bar c$ there exist $T' _1, \ldots, T' _K \in {\bf I}_{n-1} ^c(X)$, $V' _1
, \ldots , V' _K \in {\bf I}_n ^c (X)$,  $c' _1, \ldots c' _K \in C_{n-1}^L(X)$ and $m_2 \in \N_0$ with
$$\begin{array}{lcr}
 \partial \sum V' _i = T' - \sum [c' _i], && sd^{m_2}(-b \bar c ) = b \sum c'_i
\end{array}
$$
and
$$
\im(\sum c'_i) \cup \spt(\sum V' _i) \subset B_{\epsilon'}(\spt(T')\cup \im(b\bar c)) \subset  B_{R_M}(x_M) \backslash B_{\bar R}(x_M).$$
So we have an (a priori only)  integer rectifiable $n$-current with compact support
$$S \defl T\lfloor{\{d > r\}} + \sum V' _i$$
with $\spt(S) \subset \bigcup_{i=1} ^{M-1} B_{R_i}(x_i).$ The boundary is
\beq
\partial S&=& \partial (T\lfloor{\{d> r\}}) + T' - \sum [c' _i]\\
&=& (\partial T)\lfloor{\{d > r\}}-\langle T,d,r\rangle + \langle T, d, r\rangle - [\bar c]\lfloor{\{d > r\}} - \sum [c' _i]\\
&=& [c]\lfloor{\{d>  r\}} - [\bar c]\lfloor{\{d > r\}} - \sum [c' _i]\\
&=& [sd^{m_1}(c) -\bar c - \sum c' _i]  \in {\bf I}_{n-1} ^c(X),
\eeq
so $S$ is integral. By construction,  $\im(sd^{m_1}(c) -\bar c)$ and $\im(\sum c' _i)$ are subsets of $\bigcup_{i=1}^{M-1} B_{R_i} (x_i)$ and for $\tilde{c} \defl sd^{m_2}(sd^{m_1}(c) -\bar c ) - \sum c' _i$ we have  $b(\tilde{c}) = sd^{m_1 +m_2}(bc) =0$.  

By induction for $S$ 
there are $T_1, \ldots, T_{M-1}$, $c_1, \ldots, c_{M-1}$, $V_1, \ldots, V_{M-1}$ and $m_3$ with properties \ref{T2i}) and  \ref{T2ii}) of Theorem~\ref{T2} and (\ref{T2iii'}). 
Set $T_{M} \defl T - S $, then
\beq
\partial T_{M}&=& \partial T - \partial S = [c]- [\tilde{c}]\\
&=&[sd^{m_1+m_2+m_3} (c) - sd^{m_3} (sd^{m_2}(sd^{m_1}(c) -\bar c ) - \sum c' _i)]\\
&=&[sd^{m_2+m_3} (\bar c) +sd^{m_3} (\sum c' _i) ].
\eeq
With $c'' \defl sd^{m_2+m_3} (\bar c) +sd^{m_3} (\sum c' _i)$, we have $\im(c'') \subset B_{R_M}(x_M)$ and  $bc''= sd^{m_3}(sd^{m_2}(b \bar c) + b \sum c' _i)= 0$. So by the cone inequalities there exist fillings $c_{M} \in C^L_n(X)$ of $c''$ and $V_{M} \in {\bf I}_{n+1}^c(X)$ of $T-S - [c_{M}]$ with image and support
in $B_{\frac \epsilon 2} (x_M)$. Finally, 
\beq
b (\sum_{i=1}^{M} c_i)&=& sd^{m_3}(\tilde c)  + c'' \\
&=& sd^{m_1+ m_2+m_3}(c) -sd^{ m_2+m_3}(\bar c ) - sd^{m_3}(\sum c' _i)\\
&&+sd^{m_2+m_3} (\bar c) +sd^{m_3} (\sum c' _i)\\
&=& sd^m(c)
\eeq
 for $m =m_1+m_2 +m_3.$
\eprf

\subsection{Proof of Theorem~\ref{Tsing}}
The following implies Theorem~\ref{Tsing} by Proposition~\ref{contr-cone}.
\bthm     \label{Tsing2}
Let $X$ be a metric space. Then  $C_0 (X)=C_0 ^L(X)$.

Suppose the  metric space $X$ admits  local cone inequalities for $C_k(X)$ and $C_j^L(X)$ for $k=0, \dots, n$ and $j = 0, \ldots , n-1$ and let $\epsilon >0$. Given   $c \in C_n(X)$  with $bc  \in C_{n-1}^L(X)$ there exist $N \in \N$, $m \in \N_0$, $c_1, \ldots, c_N \in  C_n(X)$, $c^L_1, \dots, c^L_N \in C_n^L(X)$ and $\bar c_1, \ldots, \bar c_N \in C_{n+1} (X)$ such that
\begin{enumerate}
\item   \label{Tsing2i}
$\sum c_i = sd^m(c)$ and
$b \bar c_i = c_i- c^L _i$  
\item \label{Tsing2ii}
$\im(\bar c_i) \subset B_\epsilon (\im(c))$ and $\diam(\im(\bar c_i)) < \epsilon.$
\end{enumerate}

\ethm

The proof of this theorem is essentially the same as the proof of Theorem~\ref{T2}, when slicing is exchanged by subdivision. 
\bprf[Proof of Theorem~\ref{Tsing2}.]
Let $n=0.$ Clearly, $C_0 (X)=C_0 ^L(X)$.
Let now $n>0$ and assume that the theorem holds for $n-1$.  For $x \in X,$ let  $r_x$ denote the minimum of all  radii in the  (finitely many) cone inequalities for $x$ and $F_x$ the maximum function of all those diameter functions for $x$ 
and assume that $F_x(r) \ge r$.
Let $0 < R_x <r_x$ be such that
$$\begin{array}c
2( R_x +  F_x(2R_x + 2 F_x(2R_x))) < \epsilon.
\end{array}$$
Cover $\im(c)$ by balls with center in $\im(c)$ and of radius $R_x$ and choose a finite covering; let the centers be $x_1, \dots, x_m$, denote $R_i \defl R_{x_i}$ and $F_i \defl F_{x_i}$. 
We show by induction that there are $c_i,$ $c_i ^L$ and $\bar c_i$, $i =1, \dots , M,$ with property \ref{Tsing2i}) of Theorem~\ref{Tsing2} and the property
\begin{equation} \label{Tsing2ii'}
\im(\bar c_i) \subset B_{\frac \epsilon 2}(x_i).
\end{equation}
If $M=1$ the claim follows directly from the cone inequalities.
If $ M >1$ let  $0 < \bar R < R_M$ such that 
$\im(c) \subset B_{\bar R}(x_M) \cup \bigcup_{i=1} ^{M-1} B_{R_i} (x_i).$ Set $\alpha \defl \frac{R_M - \bar R}{4}$ and choose $r \in (\bar R + \alpha, \bar R + 3\alpha)$; note that $\bar R < r - \alpha < r +\alpha < R_M.$ Now let $m_1 \ge 0$ be such that each simplex of $sd^{m_1}(c)$ has image with diameter less than $\alpha /2.$

Let $c^+$ be the part of $sd^{m_1}(c) $ consisting of all simplices whose image has non-empty intersection with $X\backslash B_r(x_M)$ and let $c^- \defl sd^{m_1}(c) -c^+$. Let $\tilde{c}$ denote the part of the  boundary of $c^+$ that is not Lipschitz (note that is also the negative of the non-Lipschitzian part of $bc^-$) and set $\tilde{c}^L \defl bc^+ -\tilde{c}$. 

Now, $b \tilde{c} = b(bc^+- \tilde{c}^L) = -b \tilde{c}^L \in C^L _{n-2}(X)$, and $\im(\tilde{c}) \subset \bigcup_{i=1}^{M-1} B_{R_i}(x_i).$ By induction for $\tilde{c}$ with $\epsilon' \defl \alpha/2$, there are $c'_1, \ldots, c'_N \in  C_{n-1}(X)$, $c^{\prime L}_1, \dots, c^{ \prime L}_N \in C_{n-1}^L(X)$,  $\bar c_1 ', \ldots, \bar c_N' \in C_{n} (X)$ and $m_2 \in \N_0$ such that
$\sum c'_i = sd^{m_2}(\tilde{c})$, $b \bar c_i' = c_i'- c^{\prime L} _i$ and
$\im(\sum \bar c_i') \subset B_{\epsilon'} (\im(\tilde{c})) \subset B_{r+\alpha}(x_M)- B_{r-\alpha}(x_M).$

Set 
$$z \defl sd^{m_2}(c^+) - \sum \bar c_i ' \in C_n (X).$$
 Then, $\im(z) \subset \bigcup_{i=1} ^{M-1} B_{R_i}(x_i)$ and 
\beq
bz &=& sd^{m_2}(\tilde{c}^L + \tilde{c})+ \sum  c^{\prime L} _i -\sum c_i ' \\
&=& sd^{m_2}(\tilde{c}^L ) + \sum c^{\prime L} _i \in C_{n-1} ^L(X).
\eeq
By induction for $z$ 
 there are 
$c_1, \ldots, c_{M-1} \in  C_n(X)$, $c^L_1, \dots, c^L_{M-1} \in C_n^L(X)$, $\bar c_1, \ldots, \bar c_{M-1} \in C_{n+1} (X)$ and $m_3 \in N_0$ such that 
$$\sum c_i = sd^{m_3}(z),$$
$$ b \bar c_i = c_i- c^L _i \,\text{ and } \,
\im(\bar c_i) \subset B_{\frac \epsilon 2} (x_i).$$
Now, let $m \defl m_1 +m_2 +m_3$ and 
\beq
c_{M}& \defl& sd^{m}(c) -sd^{m_3}(z)\\
&=&  sd^{m}(c)-   sd^{m_2+m_3}(c^+) + sd^{m_3}(\sum \bar c_i ') \\
&=& sd^{m_2+m_3}(c^-) + sd^{m_3}( \sum \bar c_i ').
\eeq
In particular, $\im(c_{M})   \subset B_{R_M}(x_M).$ 
Denote by $c_L ^-$ the Lipschitzian part of the boundary of $c^-$ and $c_S ^- \defl bc^- - c_L^-$; recall  that $c_S^- = - \tilde{c}$ and that $\sum c_i ' = sd^{m_2}(\tilde{c})$. Then 
\beq
bc_{M} &=& sd^{m_2+m_3}(c^- _L+c_S^-) + sd^{m_3}(\sum  c_i'- \sum c^{\prime L} _i)\\
&=&  sd^{m_2+m_3}(c^- _L) -sd^{m_3}(\sum c^{\prime L} _i),
\eeq
i.e.\ $bc_{M} \in C_{n-1} ^L(X)$. The cone inequalities give now  fillings $c_{M}^L \in C_n^L(X)$ of $bc_{M}$ and $\bar c_{M} \in C_{n+1}(X)$ of $c_{M}-c_{M}^L$ such that (\ref{Tsing2ii'}) 
 is satisfied. Finally,
\beq
\sum_{i=1} ^{M} c_i &=& sd^{m_3}(z) +sd^m(c) -sd^{m_3}(z)=sd^m(c).
\eeq
\eprf

\subsection{Proof of Corollary~\ref{C4} and~\ref{C5}}
Assuming Theorem~\ref{T1} and~\ref{Tsing}, the proofs of Corollary~\ref{C4} and~\ref{C5} are the same up to minor changes. Therefore we only give the proof of Corollary~\ref{C4}.
\bprf[Proof of Corollary \ref{C4}]
We show that the chain homomorphism $[\,]{:}\, C_* ^L (X) \to {\bf I}_* ^c (X)$ induces an isomorphism $[\,]{:}\, H^ {L} _n(X, A) \to H^{IC} _n (X,A)$ for all $n \geq 0$. Note that $[\,]$ sends ${\mathcal Z}^{L} _n (X,A) $ to ${\mathcal Z}^{IC} _n (X,A)$ and  ${\mathcal B}^{L} _n (X,A) $ to $ {\mathcal B}^{IC} _n (X,A)$, so $[\,]$ induces a homomorphism from $H^ {L} _n(X, A)$ to $H^{IC} _n (X,A)$.

Let $T \in {\mathcal Z}^{IC} _n (X,A)$, as $\partial T \in {\bf I}^c _{n-1}(A)$ and $\partial(\partial T) = [0]$ there exists $c \in C_{n-1}^L(A)$ with $bc =sd^m(0) =0$ and $V \in {\bf I}^c _n(A)$ with $\partial V = \partial T - [c].$ Now, $T-V$ has boundary $\partial (T-V) = [c]$, i.e.\ there exists $ \bar c \in C_n ^L(X)$ with $b\bar c = sd^m (c)$ (thus, $\bar c \in {\mathcal Z} ^{L} _n (X,A)$) and there exists $\bar V \in {\bf I}^c _{n+1}(X)$ with $\partial \bar V = (T -V) - [\bar c]$; hence $V + \partial \bar V \in {\mathcal B}^{IC} _n (X,A)$ and 
$$[\bar c] + V + \partial \bar V = T,$$
showing the surjectivity of the homomorphism.

Let $c, \, \bar c \in {\mathcal Z}^{L} _n (X,A)$ with 
$[c] + {\mathcal B}^{IC} _n (X,A) = [\bar c] + {\mathcal B}^{IC} _n (X,A),$
 i.e.\ there exists $ R + \partial S \in {\mathcal B}^{IC} _n (X,A){:}$
$$ [c] + R + \partial S = [\bar c].$$
 Then $\partial R =  [b (\bar c -c)]$ and $b(\bar c -c ) \in C_{n-1}^L(A)$, so there exists $c_1 \in C_n^L(A), \, V \in {\bf I}^c _{n+1}(A)$ with $\partial V = R -[c_1]$ and $b c_1 = sd^{m_1}(b(\bar c - c))$. Now, by  (\ref{sd}) we get
$$b(D_{m_1,X} (\bar c -c) ) +D_{m_1,X}(b(\bar c -c))=  sd^{m_1}(\bar c - c)-\bar c + c.$$
Set $c_2 \defl D_{m_1,X} (\bar c -c) \in C_{n+1}^L(X)$; note that by naturality $D_{m_1,X}(b(\bar c -c)) =D_{m_1,A} (b(\bar c - c)) \defr c_3 \in C_n ^L(A)$.

On the other hand,
$$\partial (S+V) = [\bar c - c] -[c_1] =[sd^{m_1}(\bar c - c) -c_1].$$
Set $c_4 \defl sd^{m_1}(\bar c - c) -c_1 \in C_n ^L(X)$; then
 $bc_4 = 0$. 
So there is  a filling $c_5 \in C_{n+1}^L(X)$ with $bc_5 = sd^{m_2}(c_4).$ Thus, $ c_6\defl D_X  (c_4) \in C_{n+1}^L (X)$ has 
$$b c_6 = sd^{m_2}(c_4) - c_4 = bc_5-sd^{m_1}(\bar c - c) +c_1.$$
Together,
$$c+ c_1 -c_3 + b(c_5-c_2-c_6) = \bar c$$
with $c_1 -c_3 \in C_n ^L(A)$ and $c_5-c_2-c_6\in C_{n+1}^L(X)$, therefore the homomorphism is injective.

\eprf

\section{Proof of Theorem 1.7} \label{proofHawaiian}

The goal of this section is to show that the maximal divisible subgroup of $H_1^{IC}(\mathbb{H})$ is trivial. In \autoref{MAXIMAL_DIV_SG_NON_TRIVIAL} we will show that the maximal divisible subgroups of $H_1^L(\mathbb{H})$ and $H_1(\mathbb{H})$ are non-trivial. This implies that $H_1^{IC}(\mathbb{H})$ is not isomorphic to either one of these groups.

Recall that the Hawaiian Earring $\mathbb{H} \subset \mathbb{R}^2$ is the countable union of the circles
\[
L_n=\{x \in \mathbb{R}^2; \left\| x-(1/n,0) \right\|=1/n\}\rlap{,}
\]
with metric given by
\[
 d(x,y) \defl \left\{ \begin{matrix} \left\|x-y\right\|, & \text{ if } \exists n \in \N: \, x,y \in L_n,\;\\
                       \left\|x\right\| + \left\|y\right\|, & \text{ otherwise}
                      \end{matrix} \right.\rlap{.}
\]

\subsection{Metric Currents and the Hawaiian Earring}\label{HAWAIIAN_EARRING}

In the proof we use the fact that the first homology group of the complex of integral currents on $S^1$ is isomorphic to $\mathbb{Z}$. This follows from Corollary~\ref{C4} and~\ref{C5} since $H_1(S^1 ) \cong \Z$.

By definition, we have $\mathbf{I}_2(X)=0$ for any $\mathcal{H}^2$ null set $X$. It follows that $H_1^{IC}(\mathbb{H})$ and $H_1^{IC}(L_n)$ are simply the kernels of the maps $\partial:\mathbf{I}_1(\mathbb{H})\rightarrow \mathbf{I}_0(\mathbb{H})$ and $\partial:\mathbf{I}_1(L_n)\rightarrow \mathbf{I}_0(L_n)$ respectively. Thus, showing that an element of $H_1^{IC}(\mathbb{H})$ is zero is equivalent to showing that the integral current representing it is zero.

\begin{prop}
The maximal divisible subgroup of $H_1^{IC}(\mathbb{H})$ is trivial.
\end{prop}

\begin{proof}
Let $T \in \mathbf{I}_1 ^c(\mathbb{H})$ be an element of the maximal divisible subgroup of $H_1^{IC}(\mathbb{H})$. Let $n\in\mathbb{N}$. We write $p_n:\mathbb{H}\rightarrow L_n$ for the map which sends $x \in \mathbb{H}$ to $(0,0)$ if $x\notin L_n$ and to itself otherwise. We denote the inclusion $L_n\rightarrow \mathbb{H}$ by $i_n$. The above remark shows that $H_1^{IC}(L_n)\cong \mathbb{Z}$. We claim that $H_1^{IC}(p_n)(T)=(p_n)_{\#} (T)$ is zero. To see this, let $k\in \mathbb{Z}$ be an arbitrary integer. By assumption there exists an element $T^\prime \in H_1^{IC}(\mathbb{H})$ such that $T=k \cdot T^\prime$. It follows that $(p_n)_{\#}(T)=k\cdot (p_n)_{\#} (T^\prime) \in \mathbb{Z}$, i.e.\ every integer divides $(p_n)_{\#} (T)$, which shows that $(p_n)_{\#} (T)=0$. We find in particular that $(i_n p_n)_{\#}(T)=0$. 

On the other hand, we have $(i_n p_n)_{\#}(T)=T \lfloor L_n$. This follows since both currents have support $L_n$ and since for any function $f:\mathbb{H}\rightarrow \mathbb{R}$, the restriction of $f$ to $L_n$ agrees with the restriction of $f\circ i_n \circ p_n$ to $L_n$. These two facts imply the desired equality (by \cite[Lemma~3.2 and Theorem 4.4]{Langlocal}). 

Together this implies that
\[
\| T\|(L_n)=M(T \lfloor L_n)=M\bigl((i_n p_n)_{\#}(T)\bigr)=M(0)=0,
\]
and therefore that
\[
M(T)=\| T\|(\mathbb{H}) \leq \sum_{n=1}^\infty \|T\|(L_n)=0,
\]
which follows by countable subadditivity of $\| T\|$. We find that $T=0$, i.e.\ that the maximal divisible subgroup of $H_1^{IC}(\mathbb{H})$ is indeed trivial.
\end{proof}

\subsection{The Maximal Divisible Subgroup of \texorpdfstring{$H_1^L(\mathbb{H})$}{Lipschitz Homology}} \label{MAXIMAL_DIV_SG_NON_TRIVIAL}

\subsubsection{Overview}

In~\cite{KATSUYA_2} it is shown that the maximal divisible subgroup of $H_1(\mathbb{H})$ is non-trivial. In order to construct a nontrivial element of the maximal divisible subgroup of $H_1^L(\mathbb{H})$ we follow the construction given in the proof of Theorem~4.14 in \cite{KATSUYA_1}. Since we do not have a concise description of the Lipschitz homotopy group $\pi_1^L(\mathbb{H})$ we have to translate the algebraic definition given in loc.\ cit.\ to an explicit construction of certain Lipschitz maps $\sigma_n:[0,2\lambda(n)]\rightarrow \mathbb{H}$, $n\in\mathbb{N}$ (where $\lambda(n) \leq 1$ is a real number). The $\sigma_n$ have the following properties:
\begin{enumerate}
\renewcommand{\labelenumi}{\roman{enumi})}
\item
For $n\geq 2$, $[\sigma_{n-1}]=n\cdot [\sigma_n]$ in $H_1^L(\mathbb{H})$ and
\item
The element $[\sigma_1] \in H_1^L(\mathbb{H})$ maps to a nonzero element under the homomorphism $H_1^L(\mathbb{H})\rightarrow H_1(\mathbb{H})$.
\end{enumerate}
The idea behind the construction is the following. We first choose a sequence of maps $c_n:[0,\lambda(n)]\rightarrow \mathbb{H}$, $n\in \mathbb{N}$ which represent commutators of certain standard loops in $\pi_1(\mathbb{H})$. We construct the maps $\sigma_n$ in such a way that for $n\geq 2$ the equation
\[
\sigma_{n-1}=c_{n-1}\cdot \underbrace{\sigma_{n} \cdot \ldots \cdot \sigma_{n}}_{n \text{ times}}
\]
holds. This could be depicted as follows:

\begin{figure}[ht]
\includegraphics[width=300pt]{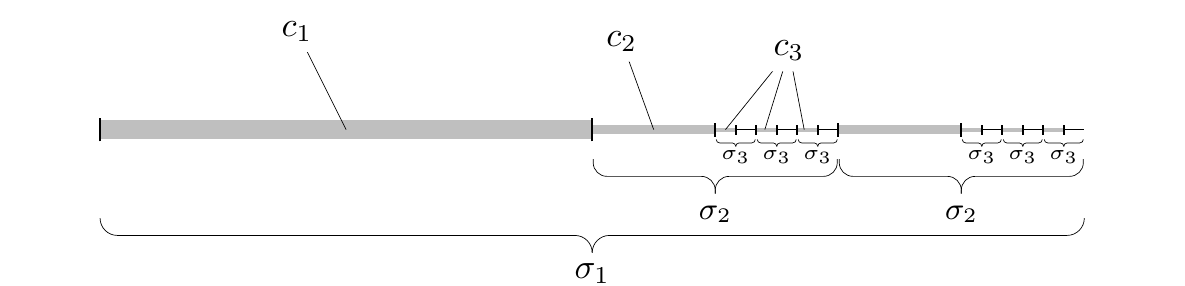}
\end{figure} 
\vspace{.5cm}
Condition i) is satisfied since $c_n$ is Lipschitz homotopic to a constant map. The commutators $c_n$ are inserted to ensure that the element $[\sigma_1]\in H_1^L(\mathbb{H})$ does not vanish, which follows from the stronger fact that its image under the comparison map $H_1^L(\mathbb{H})\rightarrow H_1(\mathbb{H})$ is a non-zero element. This is equivalent to proving that the corresponding element $[\sigma_1]_{\pi_1} \in \pi_1(\mathbb{H})$ does not lie in the commutator subgroup. We prove this in \hyperref[COMM_SUB]{Proposition~\ref{COMM_SUB}} by reducing the problem to a question about commutator subgroups of free groups.

\subsubsection{Preliminaries}
Let $S$ be the set of finite sequences of (non-zero) natural numbers, i.e.\ of maps $\{1,\ldots,n\}\rightarrow \mathbb{N}$ for $n\in \mathbb{N}$. We write $s=\langle s_1,\ldots,s_n \rangle$ for the sequence with $s(i)=s_i$. We call $n$ the \emph{length} of $s$ and denote it by $\ell(s)$. For $k,m\in \mathbb{N}$ let $k\cdot \langle m \rangle$ denote the sequence $s\in S$ of length $k$ with $s(i)=m$, $1\leq i\leq k$.
The \emph{concatenation} of sequences $+:S\times S\rightarrow S$ is given by the map which sends $(s,s^\prime) \in S\times S$ to the sequence $s^{\prime\prime}$ given by $s^{\prime\prime}(i)=s(i)$, $1\leq i \leq \ell(s)$ and $s^{\prime\prime}(j)=s^\prime(j-\ell(s))$ for $\ell(s)+1\leq j \leq \ell(s)+\ell(s^\prime)$.

Let $B=\{s\in S; s(i)\leq i \text{ and } \ell(s)>0\}$. We write $S_n$ for the subset of $S$ consisting of all sequences of length $n$ and $B_n$ for $B\cap S_n$.
There is a linear order relation $\preceq$ on $S$ such that for $s,t\in S$, $s\preceq t$ if and only if one of the following holds:
\begin{enumerate}
\renewcommand{\labelenumi}{\roman{enumi})}
\item
there exists $k \leq \min(\ell(s),\ell(t))$ such that $s(k)\neq t(k)$ and $s(j) < t(j)$ for the minimal $j$ with $s(j)\neq t(j)$
\item
$\ell(s)\leq \ell(t)$ and $s(j)=t(j)$ for all $1 \leq j\leq \ell(s)$.
\end{enumerate}
The ordered set $B$ can be embedded in the unit interval $[0,1]$. For $n\in \mathbb{N}$ let
\[
\lambda(n)=\frac{1}{2^n n!}
\]
and for $s \in B$ let
\[
\tau(s)=\sum_{t\prec s}{\lambda\bigl(\ell(t)\bigr)},
\]
where the sum is taken over any enumeration $(t_n)_{n\in \mathbb{N}}$ of $\{t\in B; t\prec s\}$ if this set is infinite. Note that the sum does not depend on the chosen enumeration and that for all $s\in B$, $\tau(s)\leq 1$. This follows since the cardinality of $B_n$ is $n!$. Indeed, let $m\in \mathbb{N}$ and write $N(m)$ for the maximum of $\{\ell(t_1),\ldots, \ell(t_m)\}$. Then the equation
\[
\sum_{k=1}^m \lambda(\ell(t_k))\leq \sum_{n=1}^{N(m)} \Bigr( \sum_{t \in B_n\text{, } t \prec s} \lambda(n) \Bigl) \leq \sum_{n=1}^\infty \lambda(n) n! = \sum_{n=1}^\infty \frac{1}{2^n}=1
\]
holds, which shows that the series is absolutely convergent. Using a suitable enumeration we find that for $s\preceq s^\prime$ the useful formula
\begin{equation}\label{differenceFormula}
\tau(s^\prime)-\tau(s)=\sum_{s\preceq t \prec s^\prime} \lambda\bigl(\ell(t)\bigr)
\end{equation}
holds.

The following lemma summarizes the basic properties of the map $\tau:B\rightarrow [0,1]$.

\begin{lemma} \label{ORDER_LEMMA}
Let $n,m\in \mathbb{N}$, $s, s^\prime \in B_n$, and let $s^{\prime\prime} \in S$ be such that $s+s^{\prime\prime}$ and $s^\prime+s^{\prime\prime}$ belong to $B$.
\begin{enumerate}[i)]
\item\label{ORDER_LEMMA_I}
The map $\tau$ is strictly order preserving. 
\item\label{ORDER_LEMMA_II}
If $s\preceq t \preceq s+\langle 1 \rangle$, then $t\in \{s,s+\langle 1 \rangle \}$. Moreover $\tau(s+\langle 1 \rangle)-\tau(s)=\lambda(n)$.
\item\label{ORDER_LEMMA_III}
Let $1\leq m \leq n$. Then $s+\langle m \rangle \preceq t \prec s+\langle m+1\rangle$ if and only if $t(i)=s(i)$ for $1\leq i \leq n$ and $t(n+1)=m$.
\item\label{ORDER_LEMMA_IV}
For $1\leq m \leq n$, $\tau(s+\langle m+1 \rangle) - \tau(s+\langle m \rangle) =2 \lambda(n+1)$.
\item\label{ORDER_LEMMA_V}
The equality $\tau(s+s^{\prime\prime})-\tau(s^\prime+s^{\prime\prime})=\tau(s)-\tau(s^\prime)$ holds. 
\end{enumerate}
\end{lemma} 

\begin{proof}
Assume $s \prec \widetilde{s}$. Then we have by definition
\[
\begin{split}
\tau(\widetilde{s})=\sum_{t \prec \widetilde{s}} \lambda\bigl(\ell(t)\bigr)&=\sum_{t \prec s} \lambda\bigl(\ell(t)\bigr) + \sum_{s\preceq t \prec \widetilde{s}} \lambda\bigl(\ell(t)\bigr) \\
&\geq \sum_{t \prec s} \lambda\bigl(\ell(t)\bigr) + \lambda\bigl(\ell(s)\bigr) > \sum_{t \prec s} \lambda\bigl(\ell(t)\bigr)=\tau(s)
\end{split}
\]
since $\lambda(\ell(s))>0$. This shows that $\tau$ is indeed strictly order preserving.

Note that the second part of ii) follows immediately from the first and equation~\eqref{differenceFormula}. Now let $s\preceq t \preceq s+\langle 1\rangle$ and let $1\leq j\leq n$. Then $s(j)<t(j)$ would imply that $s+\langle 1\rangle \prec t$, but this contradicts the hypothesis $t\preceq s+\langle 1\rangle$. Conversely, $s(j)>t(j)$ contradicts $s\preceq t$, so we find that $\ell(t)\geq n$ and $t(j)=s(j)$ for $1\leq j \leq n$. If $\ell(t)=n$ we have $t=s$, so we can assume that $\ell(t)>n$ and hence that $t(n+1)\geq 1$. If $t(n+1)$ were greater than $1$ we would have $s+\langle 1\rangle \prec t$, so we can conclude that $t(n+1)=1$. But $\ell(t)>n+1$ then implies $s+\langle 1\rangle \prec t$ (by case ii) of the definition of $\prec$) which is again a contradiction. So we finally find that $t= s+\langle 1\rangle$, which establishes the second claim.

To see iii) we have to consider two cases. First note that if $t=s+\langle m \rangle$, the conclusion holds trivially, i.e. we can assume that $s+\langle m \rangle \prec t$. Then by definition we have to consider the two cases
\begin{enumerate}[a)]
\item 
there exists $k \leq \min(n+1,\ell(t))$ such that $s+\langle m \rangle(k)\neq t(k)$ and $s+\langle m\rangle (j) < t(j)$ for the minimal $j$ with $s+\langle m\rangle (j)\neq t(j)$, or
\item
$n+1\leq \ell(t)$ and $s+\langle m \rangle (j)=t(j)$ for all $1 \leq j\leq n+1$.
\end{enumerate}
In case b) we are obviously finished, so we can assume that we are in the situation of case a). If $j \leq n$ we have $s+\langle m+1\rangle \prec t$ which contradicts our assumption. It follows that $s(i)=t(i)$ for $1\leq i \leq n$ and that $m < t(n+1)$. But this shows that $s+\langle m+1\rangle \preceq t$. Since we also have $t \prec s+\langle m+1\rangle$ we find that case a) cannot occur under the assumptions of iii).

The proof of iv) follows directly from iii) and equation~\eqref{differenceFormula}. We have
\[
\begin{split}
\tau(s+\langle m+1 \rangle) - \tau(s+\langle m \rangle)&=\sum_{t \prec s+\langle m+1 \rangle} \lambda\bigl(\ell(t)\bigr)-\sum_{t \prec s+\langle m \rangle } \lambda\bigl(\ell(t)\bigr)\\
&=\sum_{s+\langle m \rangle \preceq t \prec s+\langle m+1 \rangle} \lambda\bigl(\ell(t)\bigr)\\
&=\lambda\bigl(\ell(s+\langle m\rangle)\bigr) + \sum_{k>n+1}\Bigl(\sum_{t\in C_k} \lambda(k) \Bigr)\\
\end{split}
\]
where $C_k=\{t \in B_k; t(j)=s(j) \text{ for } 1\leq j \leq n \text{ and } t(n+1)=m\}$. Since the cardinality of $C_k$ is $\frac{k!}{(n+1)!}$ we find that
\[
\begin{split}
\tau(s+\langle m+1 \rangle) - \tau(s+\langle m \rangle)&=\lambda\bigl(\ell(s+\langle m\rangle)\bigr) + \sum_{k>n+1} \frac{k!}{(n+1)!}\lambda(k)\\
&=\lambda(n+1)+\frac{1}{(n+1)!} \sum_{k>n+1} \frac{k!}{2^k k!}\\
&=\lambda(n+1)+\frac{1}{2^{n+1}(n+1)!} \sum_{k=1}^{\infty} \frac{1}{2^k}=2\lambda(n+1).\\
\end{split}
\]

It remains to show v). There are natural numbers $m_1,\ldots, m_d$ such that $s^{\prime\prime}=\langle m_1\rangle +\ldots + \langle m_d \rangle$. By induction we can reduce the problem to the case $s^{\prime\prime}=\langle m\rangle$. Using iv) we can further reduce this to the case $m=1$. Indeed, for $m>1$ we have
\[
\begin{split}
\tau(s+\langle m\rangle)&=\tau(s+\langle m\rangle)-\tau(s+\langle m-1\rangle)+ \tau(s+\langle m-1\rangle)\\
&=\tau(s+\langle m-1 \rangle) + 2\lambda(n+1)=\ldots =\tau(s+\langle 1\rangle)+2(m-1)\lambda(n+1)
\end{split}
\]
which shows that
\[
\tau(s+\langle m\rangle)-\tau(s^\prime+\langle m\rangle)=\tau(s+\langle 1\rangle)-\tau(s^\prime+\langle 1\rangle).
\]
With ii) and equation~\eqref{differenceFormula} we find that $\tau(s+\langle 1 \rangle)=\tau(s)+\lambda(n+1)$ and therefore that the equality
\[
\tau(s+s^{\prime\prime})-\tau(s^\prime+s^{\prime\prime})=\tau(s+\langle 1\rangle)-\tau(s^\prime+\langle 1\rangle)=\tau(s)-\tau(s^\prime)
\]
holds, as claimed.
\end{proof}

\begin{lemma}\label{DENSITY_LEMMA}
The set
\[
I^\prime=\bigcup_{s\in B} [\tau(s),\tau(s+\langle 1 \rangle)] \subseteq [0,1]
\]
is dense in $[0,1]$.
\end{lemma}

\begin{proof}
Let $x\in [0,1]$ and assume that for all $s\in B$, $x \notin [\tau(s),\tau(s+\langle 1 \rangle)]$. We construct a sequence $(s_i)_{i\in\mathbb{N}}$ with $\tau(s_i)<x$ as follows. Let $s_1=\langle 1 \rangle$ (which trivially satisfies $\tau(s_1)=0 < x$), and if $s_k$ is defined for all $k<i$, let $s_i$ be the maximal element of $\{s\in B_i; \tau(s)< x\}$. This set is not empty since it contains $s_{i-1}+\langle 1 \rangle$. Indeed, by induction hypothesis we have $\tau(s_{i-1})<x$ and $x\notin I^\prime$ implies in particular $x \notin [\tau(s_{i-1}),\tau(s_{i-1}+\langle 1 \rangle)]$, so $\tau(s_{i-1}+\langle 1 \rangle)< x$. Moreover, for any $i\in \mathbb{N}$ there exists $m\in \mathbb{N}$ such that
\begin{equation}\label{S_I_FORMULA}
s_{i+1}=s_i+\langle m \rangle.
\end{equation}
This follows since $s_i \prec s_i+\langle 1 \rangle \preceq s_{i+1}$ and therefore one of the following must hold:
\begin{enumerate}[a)]
\item
there is a $k\leq i$ with $s_i(k)\neq s_{i+1}(k)$ and $s_i(j) < s_{i+1}(j)$ for the minimal $j$ with $s_i(j)\neq s_{i+1}(j)$, or
\item 
the sequence $s_{i+1}$ extends $s_i$, i.e. $s_{i+1}=s_i+\langle m \rangle$.
\end{enumerate}
But case a) would contradict the maximality of $s_i$: $s^\prime=\langle s_{i+1}(1),\ldots , s_{i+1}(i)\rangle$ is an element of $B_i$, the equation
\[
\tau(s^\prime)\leq \tau\bigl(s^\prime+\langle s_{i+1}(i+1)\rangle\bigr)=\tau(s_{i+1}) < x
\]
holds, and $s^\prime \prec s_i$ since $s_i(j)<s_{i+1}(j)=s^\prime(j)$. 

This implies in particular that $(\tau(s_i))_{i\in\mathbb{N}}$ is a strictly increasing sequence and that $\sup \{\tau(s_i);i\in\mathbb{N}\}=\lim_{i\rightarrow \infty} \tau(s_i) \leq x$. We will show that this is an equality and therefore that $x$ lies in the closure of $I^\prime$.

First note that the set
\[
C=\{i\in \mathbb{N};s_i(i)<i\}
\]
is either empty or infinite. If $C$ were finite and non-empty, it would contain a maximal element $n\in \mathbb{N}$. We write $s^\prime$ for the sequence $\langle s_n(1), s_n(2), \ldots, s_n(n-1) \rangle$ and let $m=s_n(n)$. Since $n \in C$ we have $m<n$. By maximality of $n$ and with equation~\eqref{S_I_FORMULA} it follows that for $i>n$ we must have $s_i=s_n +\langle n+1, n+2, \ldots, i \rangle$. Now let $t\prec s^\prime +\langle m+1 \rangle$. Then we have either $t\prec s_n$ or $t(j)=s^\prime(j)$ for $1\leq j \leq n-1$ and $t(n)=m$ (by \hyperref[ORDER_LEMMA_III]{Lemma~\ref{ORDER_LEMMA}~\ref{ORDER_LEMMA_III})}). But in the latter case we must have $t\preceq s_{\ell(t)}$ since, for any $i>n$, $s_i$ is the maximal element of $B_i$ satisfying the given condition. So in either case there is an $i\in \mathbb{N}$ with $t\prec s_i$. It follows that
\[
\begin{split}
\tau(s^\prime + \langle m+1\rangle) &=  \sum_{t\prec s^\prime+\langle m+1 \rangle} \lambda\bigl(\ell(t)\bigr) \leq \sup \biggl\lbrace\sum_{t^\prime \preceq t} \lambda\bigl(\ell(t^\prime)\bigr); t\prec s^\prime+\langle m+1 \rangle \biggr\rbrace\\
&\leq \sup \biggl\lbrace\sum_{t^\prime \prec s_k} \lambda\bigl(\ell(t^\prime)\bigr); k \in \mathbb{N} \biggr\rbrace=\lim_{i\rightarrow \infty} \tau(s_i) \leq x\\
\end{split}
\]
and hence, because $x \in I^{\prime}$, that $\tau(s^\prime + \langle m+1\rangle)<x$. But this contradicts the maximality of $s_n$ in $\{s\in B_n; \tau(s)< x\}$, so $C$ must be infinite.

We first consider the case $C=\varnothing$. By construction we must have $s_i=\langle 1, 2, \ldots, i \rangle$ and it follows that
\[
\tau(s_i) \geq \sum_{k=1}^{i-1}\bigl( \sum_{t\in B_k} \lambda(k) \bigr)=\sum_{k=1}^{i-1} \frac{1}{2^k}=1-\frac{1}{2^{i-1}}\\,
\]
which shows that $\lim_{i\rightarrow \infty} \tau(s_i)=1$. But this implies that $x=1$ and that $x$ is the limit of $(\tau(s_i))_{i \in \mathbb{N}}$. 

It remains to show that the same holds if $C$ is infinite. We prove this by contradiction, so assume now that $\lim_{i\rightarrow \infty} \tau(s_i)<x$. Choose $n\in \mathbb{N}$ such that $2\lambda(n) < \varepsilon$, where $\varepsilon=x- \lim_{i\rightarrow \infty} \tau(s_i)$. Since the set $C$ is infinite it follows that there exists $N>n$ such that $s_N \in C$. But now \hyperref[ORDER_LEMMA_IV]{Lemma~\ref{ORDER_LEMMA}~\ref{ORDER_LEMMA_IV})} implies that we have
\[
\tau(\langle s_N(1), s_N(2), \ldots s_N(N-1) \rangle +\langle s_N(N)+1\rangle)=\tau(s_N)+2\lambda(N) < \tau(s_N)+\varepsilon \leq x,
\]
which contradicts the maximality of $s_N$ in $\{s\in B_N; \tau(s)< x\}$. So we have indeed found a sequence $(\tau(s_i))_{i\in \mathbb{N}}$ in $I^\prime$ with $\lim_{i\rightarrow \infty} \tau(s_i)=x$.
\end{proof}

This fact enables us to specify Lipschitz functions from $[0,1]$ to a complete metric space $X$ (e.g. $\mathbb{H}$) by specifying their restriction to $I^\prime$. This is quite simple since for $s \prec s^\prime \in B$, the intersection of the two intervals $[\tau(s),\tau(s+\langle 1 \rangle)]\cap [\tau(s^\prime),\tau(s^\prime+\langle 1 \rangle)]$ is either empty or the singleton $\{\tau(s^\prime)\}$. So if we have a $1$-Lipschitz function
\[
\varphi_s:[\tau(s),\tau(s+\langle 1 \rangle)]\rightarrow \mathbb{H}
\] 
for each $s \in B$ with $\varphi_s(\tau(s))=\varphi_s(\tau(s+\langle 1 \rangle))=x_0$, there is a unique $1$-Lipschitz function $\varphi:[0,1]\rightarrow \mathbb X$ with
\begin{equation}
\varphi\vert_{[\tau(s),\tau(s+\langle 1\rangle)]}=\varphi_s.
\end{equation}
We use this construction to define $\sigma_1:[0,1]\rightarrow \mathbb{H}$.

The following construction is useful for the computation of certain finite sums in $H_1^L(X)$. Given two Lipschitz maps $\sigma:[0,a]\rightarrow X$ and $\sigma^\prime:[0,b]\rightarrow X$, $a,b\in\mathbb{R}$ with $\sigma(a)=\sigma^\prime(0)$, their \emph{concatenation} $\sigma\cdot \sigma^\prime:[0,a+b]\rightarrow X$ is given by
\[
\sigma\cdot \sigma^\prime(t)=\biggl\{ 
\begin{array}{ll}
\sigma(t) & \text{if } t\leq a\\
\sigma^\prime(t-a) & \text{if } a\leq t \leq  b\\
\end{array}.
\]
If $\sigma$ and $\sigma^\prime$ represent elements of $H_1^L(X)$, so does $\sigma\cdot \sigma^\prime$ and moreover the equation
\[
[\sigma\cdot \sigma^\prime]=[\sigma]+[\sigma^\prime]
\]
holds. This follows as in the continuous case.

\subsection{Construction of the \texorpdfstring{$\sigma_n$}{Sigma's}}

We write $\varphi_n$ for the Lipschitz function $[0,1]\rightarrow \mathbb{H}$ which traverses the $n$-th loop $L_n$ of $\mathbb{H}$ with constant speed. Explicitly, for $t\in [0,1]$ the equation
\begin{equation}\label{PHI_N}
\varphi_n(t)=\Biggl(-\frac{\cos(2\pi t)}{n}+\frac{1}{n},\frac{\sin(2\pi t)}{n}\Biggr)
\end{equation}
holds. The Lipschitz constant of $\varphi_n$ is bounded by $\mu_n=2\pi/n$. Choose a sequence $(n_k)_{k\in\mathbb{N}}$ such that
\[
4\mu_{n_k} \leq \lambda(k)
\]
and $n_{k+1}>n_k+1$. We write 
\begin{equation}\label{C_K}
c_k=\varphi_{n_k}\cdot \varphi_{n_k+1} \cdot \varphi_{n_k}^{-1} \cdot \varphi_{n_k+1}^{-1} \circ \psi:[0,\lambda(k)]\rightarrow \mathbb{H}
\end{equation}
for the composition of (a representative of) the commutator of $\varphi_{n_k}$ and $\varphi_{n_k+1}$ with the  reparametrization $\psi:[0,\lambda(k)]\rightarrow [0,4]$ which sends $t$ to $4t/\lambda(k)$. By choice of the sequence $(n_k)_{k\in\mathbb{N}}$ we find that $c_k$ is a $1$-Lipschitz function. For $s\in B$ and $t\in [\tau(s),\tau(s+\langle 1 \rangle)]$ let $(\sigma_1)_s(t)=c_{\ell(s)}(t-\tau(s))$ and note that $(\sigma_1)_s(\tau(s))=(\sigma_1)_s(\tau(s+\langle 1 \rangle))=(0,0)$. By the comment succeeding \hyperref[DENSITY_LEMMA]{Lemma~\ref{DENSITY_LEMMA}} there is a unique $1$-Lipschitz function $\sigma_1:[0,1]\rightarrow \mathbb{H}$ such that
\begin{equation}\label{SIGMA_1}
\sigma_1\vert_{[\tau(s),\tau(s+\langle 1 \rangle)]}(t)=c_{\ell(s)}(t-\tau(s)).
\end{equation}
holds. For $n>1$ Let $\sigma_n:[0,2\lambda(n)]\rightarrow \mathbb{H}$ be the function which sends $t \in [0,2 \lambda(n)]$ to
\begin{equation}\label{SIGMA_N}
\sigma_n(t)=\sigma_1\vert_{[\tau(n\cdot\langle 1 \rangle),\tau((n-1)\cdot\langle 1 \rangle+\langle 2\rangle )]}(t+\tau(n\cdot\langle 1 \rangle)).
\end{equation}

\begin{prop}
For $n>1$ the equation
\[
\sigma_{n-1}=c_{n-1} \cdot \underbrace{\sigma_n \cdot \ldots \cdot \sigma_n}_{n \text{ \normalfont{times}}}
\]
holds, and consequently $[\sigma_{n-1}]=n\cdot [\sigma_n]$ in $H^L_1(\mathbb{H})$. 
\end{prop}

\begin{proof}
Let $a=\tau\bigl((n-1)\cdot\langle 1 \rangle\bigr)$, $t \in [0,2\lambda(n-1)]$ and let $x=t+a$. If $t\leq\lambda(n-1)$ or equivalently $x \in \bigl[\tau\bigr((n-1)\cdot \langle 1 \rangle\bigl),\tau\bigr(n\cdot \langle 1 \rangle\bigl)\bigr]$ we find that
\[
\begin{split}
\sigma_{n-1}(t)&=\sigma_1(x)=\sigma_1\vert_{\bigl[\tau((n-1)\cdot \langle 1 \rangle),\tau(n\cdot \langle 1 \rangle)\bigr]}(x)\\
&=c_{n-1}(x-a)=c_{n-1}(t)=c_{n-1}\cdot \underbrace{\sigma_n \cdot \ldots \cdot \sigma_n}_{n \text{ \normalfont{times}}}(t)\\
\end{split}
\]
holds. We can therefore restrict our attention to those points $t\in [\lambda(n-1),2\lambda(n-1)]$ with the property that $t+a$ lies in $I^\prime$, so we can assume that there is a sequence $s_0 \in B$ such that $x=t+a$ lies in the interior of $\bigl[\tau(s_0),\tau\bigl(s_0+\langle 1 \rangle\bigr)\bigr)$. But $t\geq \lambda(n-1)$ implies that $x\geq \tau(n\cdot\langle 1 \rangle)$, and it follows that $n\cdot \langle 1 \rangle \preceq s_0$ since $s_0$ is maximal with $\tau(s_0)\leq x$ (see \hyperref[ORDER_LEMMA_II]{Lemma~\ref{ORDER_LEMMA}~\ref{ORDER_LEMMA_II})}). On the other hand, we have
\[
\tau(s_0)<x\leq \tau\bigl((n-1)\cdot \langle 1\rangle \bigr)+2\lambda(n-1)=\tau((n-2)\cdot \langle 1\rangle +\langle 2\rangle)
\]
and therefore $(n-1)\cdot \langle 1 \rangle \preceq s_0 \prec (n-2)\cdot \langle 1\rangle +\langle 2\rangle$. By \hyperref[ORDER_LEMMA_III]{Lemma~\ref{ORDER_LEMMA}~\ref{ORDER_LEMMA_III})} it follows that $s_0=(n-1)\cdot\langle 1\rangle + s^\prime$ for some sequence $s^\prime \in S$. But $n\cdot \langle 1 \rangle \preceq (n-1) \cdot\langle 1\rangle + s^\prime$ implies that $s^\prime$ is of the form $\langle m \rangle + s$ for some $m\in\mathbb{N}$ and a (possibly empty) sequence $s\in S$. Hence $s_0=(n-1)\cdot\langle 1\rangle + \langle m \rangle + s$ and we find that
\[
\begin{split}
\sigma_{n-1}(t)&=\sigma_1(x)=c_{n+\ell(s)}\Bigl(x-\tau\bigl((n-1)\cdot \langle 1 \rangle + \langle m \rangle + s\bigr)\Bigr)\\
&=c_{n+\ell(s)}\Bigl(x-\tau\bigl((n-1)\cdot \langle 1 \rangle + \langle m \rangle + s\bigr)+\tau\bigl(n\cdot \langle 1 \rangle + s\bigr)-\tau\bigl(n\cdot \langle 1 \rangle + s\bigr)\Bigr)\\
&=c_{n+\ell(s)}\Bigl(x-\tau\bigl((n-1)\cdot \langle 1 \rangle + \langle m \rangle\bigr)+\tau\bigl(n\cdot \langle 1 \rangle\bigr)-\tau\bigl(n\cdot \langle 1 \rangle + s\bigr)\Bigr)\\
&=\sigma_1\Bigl(x-\tau\bigl((n-1)\cdot \langle 1 \rangle + \langle m \rangle\bigr)+\tau\bigl(n\cdot \langle 1 \rangle\bigr)\Bigl)\\
&=\sigma_n\Bigl(x-\tau\bigl((n-1)\cdot \langle 1 \rangle + \langle m \rangle\bigr)\Bigl)\\
&=\sigma_n\Bigl(t+\tau\bigl((n-1)\cdot\langle 1 \rangle\bigr)-\tau\bigl((n-1)\cdot \langle 1 \rangle + \langle m \rangle\bigr)\Bigl)\\
&=\sigma_n\bigl(t-\lambda(n-1)-(m-1)\cdot 2\lambda(n)\bigr)\\
\end{split}
\]
holds, where the crucial step is an application of \hyperref[ORDER_LEMMA_V]{Lemma~\ref{ORDER_LEMMA}~\ref{ORDER_LEMMA_V})} and the last equality follows from \hyperref[ORDER_LEMMA_II]{Lemma~\ref{ORDER_LEMMA}~\ref{ORDER_LEMMA_II})} and \hyperref[ORDER_LEMMA_III]{\ref{ORDER_LEMMA_III})}. On the other hand, since we have
\[
\tau\bigl((n-1)\langle 1 \rangle + \langle m \rangle \bigr) \leq x \leq \tau\bigl((n-1)\langle 1 \rangle + \langle m+1 \rangle \bigr)
\]
it follows that $\lambda(n-1)+(m-1)\cdot 2\lambda(n) \leq t \leq \lambda(n-1)+m\cdot2\lambda(n)$, again from \hyperref[ORDER_LEMMA_II]{Lemma~\ref{ORDER_LEMMA}~\ref{ORDER_LEMMA_II})} and \hyperref[ORDER_LEMMA_III]{\ref{ORDER_LEMMA_III})}. So in this case we have
\[
\begin{split}
c_{n-1} \cdot \underbrace{\sigma_n \cdot \ldots \cdot \sigma_n}_{n \text{ \normalfont{times}}}(t)&=\underbrace{\sigma_n \cdot \ldots \cdot \sigma_n}_{n \text{ \normalfont{times}}}(t-\lambda(n-1))\\
&=\underbrace{\sigma_n \cdot \ldots \cdot \sigma_n}_{n-1 \text{ \normalfont{times}}}(t-\lambda(n-1)-2\lambda(n))\\
&=\underbrace{\sigma_n \cdot \ldots \cdot \sigma_n}_{n-2 \text{ \normalfont{times}}}(t-\lambda(n-1)-2\cdot 2\lambda(n))\\
&\quad\quad\quad\quad \vdots\\
&=\underbrace{\sigma_n \cdot \ldots \cdot \sigma_n}_{\bigl(n-(m-1)\bigr) \text{ \normalfont{times}}}(t-\lambda(n-1)-(m-1)\cdot 2\lambda(n))\\
&=\sigma_n(t-\lambda(n-1)-(m-1)\cdot 2\lambda(n)),\\
\end{split}
\]
and we find that the claimed equality holds for a dense subset.
The second statement follows since we have in general $[\sigma\cdot\sigma^\prime]=[\sigma]+[\sigma^\prime]$. This implies in particular that the commutator $[c_{n-1}]$ vanishes (since $H^L_1(\mathbb{H})$ is abelian) and therefore that $[\sigma_{n-1}]=n\cdot[\sigma_n]$.
\end{proof}

\begin{cor}
The element $[\sigma_1]$ lies in the maximal divisible subgroup of $H^L_1(\mathbb{H})$.
\end{cor}

\begin{proof}
The set
\[
D=\{m \cdot [\sigma_n]; m\in\mathbb{Z},n\in \mathbb{N}\}
\]
is a divisible subgroup of $H^L_1(\mathbb{H})$. To see this, let $m,m^\prime \in\mathbb{Z}$, $n,n^\prime \in\mathbb{N}$ with $n\leq n^\prime$. Then $[\sigma_n]=\frac{n^\prime!}{n!}[\sigma_{n^\prime}]$ by the above proposition. It follows that
\[
m \cdot [\sigma_n]-m^\prime \cdot [\sigma_{n^\prime}]=\biggl(m\cdot\frac{n^\prime!}{n!}-m^\prime\biggr)[\sigma_{n^\prime}],
\]
so $D$ is indeed a subgroup. But it is also divisible, for if $l\in \mathbb{N}$ we have
\[
m \cdot [\sigma_n]=m\cdot \frac{(l\cdot n)!}{n!}[\sigma_{l\cdot n}]=l\cdot \biggl(m\cdot \frac{(l\cdot n-1)!}{(n-1)!}\biggr)[\sigma_{l\cdot n}]
\]
and $m\cdot \frac{(l\cdot n-1)!}{(n-1)!}[\sigma_{l\cdot n}] \in D$.
\end{proof}

It remains to show that $[\sigma_1]\neq 0$. We will show that the image $[\sigma_1]_S \in H_1(\mathbb{H})$ of $[\sigma_1]$ under the comparison map $H_1^L(\mathbb{H})\rightarrow H_1(\mathbb{H})$ is non-zero. Using the isomorphism
\[
H_1(\mathbb{H}) \cong \pi_1(\mathbb{H})/[\pi_1(\mathbb{H}),\pi_1(\mathbb{H})]
\]
we find that this is equivalent to showing that the loop $\sigma_1$ does not represent an element of the commutator subgroup of $\pi_1(\mathbb{H})$.

Let $p_k:\mathbb{H}\rightarrow L_{n_k}\cup L_{n_k+1}$ be the map which is the identity on $L_{n_k}\cup L_{n_k+1}$ and which sends an element $x \in \mathbb{H}\setminus L_{n_k}\cup L_{n_k+1}$ to $(0,0)$. This is a (Lipschitz) continuous map and therefore induces a group homomorphism $\pi_1(p_k):\pi_1(\mathbb{H})\rightarrow \pi_1(L_{n_k}\cup L_{n_k+1})$. Since $L_{n_k}\cup L_{n_k+1}$ is homeomorphic to $S^1\vee S^1$, $\pi_1(L_{n_k}\cup L_{n_k+1})$ is isomorphic to the free group on two generators $\langle a,b\rangle$. Equation \eqref{PHI_N} implies that $\pi_1(p_k)([\varphi_{n_k}]_{\pi_1})=[p_k\circ \varphi_{n_k}]_{\pi_1}=a$ and $\pi_1(p_k)([\varphi_{n_k+1}]_{\pi_1})=[p_k\circ \varphi_{n_k+1}]_{\pi_1}=b$ under this isomorphism. We find in particular that
\begin{equation}\label{COMMUTATOR_EQN}
\pi_1(p_k)([c_k])=aba^{-1}b^{-1}.
\end{equation}

\begin{prop}\label{COMM_SUB}
The element $[\sigma_1]_{\pi_1} \in \pi_1(\mathbb{H})$ does not lie in the commutator subgroup $[\pi_1(\mathbb{H}),\pi_1(\mathbb{H})]$.
\end{prop}

\begin{proof}
>From now on we identify $\pi_1(L_{n_k}\cup L_{n_k+1})$ with $\langle a,b\rangle$. We will first show that $\pi_1(p_k)$ maps $[\sigma_1]_{\pi_1}$ to $(aba^{-1}b^{-1})^{k!}$. By \hyperref[DENSITY_LEMMA]{Lemma~\ref{DENSITY_LEMMA}} the map $p_k\circ \sigma_1$ is uniquely determined by its restriction to the intervals $[\tau(s),\tau(s+\langle 1 \rangle)]$. By equation $\eqref{C_K}$ it follows that
\[
p_k\circ \sigma_1\vert_{[\tau(s),\tau(s+\langle 1 \rangle)]}(t)=p_k \circ c_{\ell(s)}(t-\tau(s))=\biggl\{ 
\begin{array}{lr}
c_k(t-\tau(s)) & \text{if } \ell(s)=k\\
(0,0) & \text{else}\\
\end{array}.
\]
Since there are precisely $k!$ sequences $s\in B$ with $\ell(s)=k$ we find that $p_k\circ \sigma_1$ is homotopic (rel.\ endpoints) to
\[ 
\underbrace{c_k\cdot \ldots \cdot c_k}_{k! \text{ times}}
\]
and together with equation~\eqref{COMMUTATOR_EQN} that $\pi_1(p_k)([\sigma_1]_{\pi_1})=(aba^{-1}b^{-1})^{k!}$.

Assume now that $[\sigma_1]_{\pi_1}$ does lie in the commutator subgroup, i.e.\ that $[\sigma_1]_{\pi_1}=[x_1,y_1]\cdot \ldots \cdot [x_n,y_n]$ for some $x_i,y_i\in \pi_1(\mathbb{H})$. Since $\pi_1(p_k)$ is a group homomorphism, it follows that for all $k\in\mathbb{N}$ the element $(aba^{-1}b^{-1})^{k!}$ of $\langle a, b\rangle$ can be written as a product of  $n$ elementary commutators. On the other hand, in \cite[example~2.6]{CULLER}  it was shown that $(aba^{-1}b^{-1})^{k!}$ cannot be written as a product of less than $\lfloor k!/2\rfloor + 1$ elementary commutators.
Since $k\in \mathbb{N}$ was arbitrary this is clearly a contradiction.
\end{proof}

\bibliographystyle{alpha}
\bibliography{Homology}
\end{document}